\documentclass[11pt]{article}
\usepackage{amsmath,amssymb,amsthm,mathrsfs}
\usepackage{hyperref}

\begin{document}

\renewcommand{\baselinestretch}{1.3}
\renewcommand{\arraystretch}{1.3}

\begin{center}{\bf\LARGE  Frobenius Theorem in Banach Space and Generalized Inverse Analysis Method of Operators Under Small Perturbations   }\\
\vskip 0.5cm Ma Jipu\\[10pt]
Tseng Yuanrong Functional Research Center, Harbin Normal
  University, Harbin, 150080, P. R. China\\

Department of Mathematics, Nanjing University, Nanjing, 210093,
  P. R. China

 E-mail address: jipuma@126.com

\end{center}

{\bf Abstract}\quad   Let $\Lambda$ be an open set in Banach space $E$, $M(x)$ for $x\in \Lambda $ be a subspace in $E$, and $x_0$ be a point in $\Lambda $. We consider the family $\mathcal{F}=\{M(x):\forall x\in\Lambda\}$, but the dimension of $M(x)$ can be infinite, and investigate the necessary and sufficient conditions for $\mathcal{F}$ being $c^1$ integrable at $x_0$. Without new idea and method, it is difficult to generalize the classical Frobenius theorem in Euclid space to the infinite-dimensional $M (x)$ case. We first define the co-tailed set $J (x_0, E_ *)$ of $\mathcal{F}$ at $x_0$ so that for each $x$ in $J (x_0, E_ *)$, $M (x)$ has a unique operator value coordinate $\alpha(x)$ in $B(M (x_0), E_*),$ and prove that if $\mathcal{F}$ is integrable at $x_0$, $J (x_0, E_ *)$ must contain the integrable submanifold of $\mathcal{F}$ at $x_0$. Then, we present the desired necessary and sufficient conditions, which is the Frobenius theorem in the Banach space.It is well known that the classical Frobenius theorem is an important fundamental theorem in the fields of differential topology, differential geometry, differential equations, etc. However, they are all limited to cases where all $\mbox{dim}M(x)< \infty.$ It is now possible to generalize previous studies to the case of $\mbox{dim} M(x)=\infty.$
Using the generalized inverse analysis method of operators under small perturbations, we not only prove Frobenius theorem, but also give some applications to the initial value problem of differential equations with geometric significance, global analysis and the extremum principle under the submanifold constraint in Banach space. In particular, in the field of infinite dimensional geometric and functional analysis, these studies seem to belong to new results and are still in the preliminary stage.



\vskip 0.2cm\begin{center}{\bf 1\quad Introduction and Preliminary }\end{center}\vskip 0.2cm

 Let $E$ be a Banach space, and $\Lambda$ an open set in $E$. Assign a subspace $M(x)$ in $E$ to each point $x$ of $\Lambda$ , but
$\dim M(x)$ can be infinite. We consider the family  $\mathcal{F}=\{M(x):\forall x\in\Lambda\}$, and investigate the
necessary and sufficient conditions for $\mathcal{F}$ being $c^1$ integrable at a point $x_0\in\Lambda$.We first introduce the co-tailed set $J (x_0, E_ *)$ of $\mathcal{F}$ at $x_0$ (see Definition 2.1)so that for each $x$ in $J (x_0, E_ *)$, $M(x)$ has a unique operator-valued coordinate $\alpha(x)$ in $B(M (x_0), E_*),$ satisfying $M (x) = \{e + \alpha(x)e : e \in M (x_0)\}$(see Theorem 1.5). Through $J (x_0, E_ *)$ and $\alpha(x)$, we ultimately establish Frobenius' theorem in the Banach space(see Theorem $2.2$).
For example, we apply the generalized inverse analysis method of operators under small perturbations to initial valued problem of differential equations with geometric meaning,to global analysis, and to the extremum principle under constraints of $c^1$ submanifolds in Banach space. (See Theorem $3.2$ and Theorem $4.1-4.4$.) In recent years, the rank theorem of matrices has been extended to operators and applied in nonlinear functional analysis, thus forming the generalized inverse analysis method of operators under small perturbations. At present, most of our colleagues are not familiar with this content, so we introduce Theorem $1.1-1.5$ in this section and give their proofs in the appendix.Recall that $A^+\in B(F, E)$ is said to be a generalized inverse of $A\in B(E,F)$\ provided $A^+AA^+=A^+$ and $A=AA^+A$; $A\in B(E,F)$ is said to be double-split if $R(A)$ is closed and there is a closed subspace $R^+$ in $E$, and a closed subspace $N^+$ in $F$, such that $E=N(A)\oplus R^+$ and $F=R(A)\oplus N^+.$ It is well known that $A$ has a generalized inverse $A^+\in B(F,E)$ if and only if $A$ is double-split (see [Ma.1] and [N-C]). In the sequel,let $B^+(E,F)$ denote the set of double-split operators in $B(E,F),$ and $A^+$ be a generalized inverse of $A$ in  $B^+(E,F)$.For a non zero operator $A$ in $B^+(E,F)$, write
  $V(A,A^+)=\{T\in B(E,F):\|T-A\|<\|A^+\|^{-1}\}$, $C_A(A^+,T)=I_F+(T-A)A^+$ and
$D_A(A^+,T)=I_E+A^+(T-A).$ Then we have

{\bf Theorem 1.1.}\quad {\it For $T \in V(A,A^+),$ the following conditions are  mutually equivalent:

(i) $R(T)\cap N(A^+)=\{0\}$;

(ii) $B=A^+C^{-1}_A(A^+,T)= D^{-1}_A(A^+,T)A^+) $ is the
generalized inverse of $T$, which meets both $R(B)=R(A^+)$ and $N(B)=N(A^+)$ (see [N-C]);

(iii) $R(T)\oplus N(A^+)=F;$

(iv) $N(T)\oplus R(A^+)=E$;

(v) $(I_E-A^+A)N(T)=N(A)$;

(vi) $C^{-1}_A(A^+,T)TN(A)\subset R(A)$;

(vii) $R(C^{-1}_A(A^+,T)T)\subset R(A)$.

(See [N-C], [Ma1] and [Ma.4].)}

Let $F_k=\{A\in B^+(E,F): \mbox{dim} R(A)=k\},$ $\Phi_{m,n}=\{A\in B^+(E,F): \mbox{dim} N(A)=m$  $\mbox{and} \, \, \mbox{codim} R(A)=n\},$
$\Phi_{m,\infty}=\{A\in B^+(E,F): \mbox{dim} N(A)=m \, \, \mbox{and} \, \, \mbox{codim} R(A)=\infty,\}$ and $\Phi_{\infty,n}=\{A\in B^+(E,F): \mbox{dim} N(A)=\infty \, \, \mbox{and} \, \, \mbox{codim} R(A)=n\}.$ We have

 {\bf Theorem 1.2.}\quad {\it  Assume that $A$ belongs to any of the following classes:
$$ F_k, \, \Phi_{m,n}, \, \Phi_{m,\infty}, \, \, \mbox{and} \quad \Phi_{\infty, n} ,$$
$ k=1,2,\cdots $ and $ m,n=0,1, \cdots.$ Let $T$ be an operator in $V(A^+,A)$. The condition $R(T)\cap N(A^+)=\{0\} $ holds
if and only if $T$ belongs to the same class with $A.$ (See [Ma 2] and [Ma 4].)}

Let $T_x$ be a continuous operator-valued mapping from the topological space $X$ to $B(E,F)$ at point $x_0$.

{\bf Definition 1.1.}\quad The point $x_0$ is said to be a locally fine point of $T_x$
 provided that for a generalized inverse $T^+_{x_0}$ of
 $T_{x_0},$ there exists a neighborhood $U_0$ at $x_0$,
 such that
 $$R(T_x)\cap N(T^+_{x_0})=\{0\} \,\, \mbox{for all} \,\, x\in U_0$$
where the neighborhood $U_0$ generally depends on $T^+_{x_0}.$

Let $f$ be a $c^1$ mapping from an open set $U$ in $E$ to the space $F$, and $x_0$ a point in $U.$

{\bf Definition 1.2.}\quad{\it  The locally fine point
 $x_0$ of $f'(x)$ is said to be a generalized regular point of $f$.} (See [Ma 1-4].)

Theorem $1.2$ provides generalized regular points of various $c^1$ maps $f$. (See [Ma $1$], [Ma$2$], [Abr] and [Zei].)

{\bf Theorem 1.3.}\quad{\it  The locally fine point
 $x_0$ of $T_x$ is independent of the choice of generalized inverse $T^+_{x_0}$ of $T_{x_0}.$}

{\bf Note:} local fine points and generalized regular points are two key concepts in the generalized inverse analysis method of operators under small perturbations, which bring many new theorems in nonlinear analysis and submanifolds in global analysis.(See Theorem 1.4, Theorem 3.2, Theorems 4.1-4.3, and see [Ma1], [Ma3], [Ma6], [Ma8], and [Ma9].)

{\bf Theorem 1.4.}\quad (Operator rank theorem)\quad {\it Suppose that
 $T_{x_0}$ is double-split. Then for any generalized inverse $T^{+}_{x_0}$ of
$T_{x_0}$, there exists a neighborhood $U_0$ at $x_0$ that satisfies the following two conditions: for each $x\in U_0$, there is a generalized inverse $T^+_x$ of $T_x$ and $\lim\limits_{x\rightarrow x_0}T^+_x=T^+_{x_0}$, if and only if $x_0$ is a local fine point of $T_x$.}

The famous rank theorem of matrices states that when $\Delta A$ is small enough,if $Rank (A)=Rank (A+\Delta A)$, then as $\Delta A$ approaches zero, $(A+\Delta A)^+$ approaches $A^+$, where $A$ and $\Delta A$ are both $(n\times m)$ matrices. The rank theorem guarantees the stability of computing generalized inverse matrices, therefore, it is a very important theorem in computational mathematics. (See [Ben] and [P].)Theorem 1.4 extends the rank theorem of matrices to bounded linear operators between Banach spaces, which in particular leads to some important developments in nonlinear analysis theory (see [Ma1],[Ma3], [Ma5], [Ma7] and [MA8]).

{\bf Theorem 1.5.}\quad {\it Suppose that the two closed subspaces $E_0$ and $E_1$ in $E$ have a common complement
$E_*$. Then there exists a unique operator $\alpha\in B(E_0,E_*)$, such that
$$E_1=\{e+\alpha e:\forall e\in E_0\}.$$ Further,
for any $\alpha\in B(E_0,E_*)$, $$\{e+\alpha e:\forall e\in E_0\}\oplus E_*=E .$$}

For detail proofs of Theorems $1.1-1.5,$ see the appendix in this
paper.

 \vskip 0.2cm\begin{center}{\bf 2\quad  Frobenius
 Theorem in Banach Space }\end{center}\vskip 0.2cm

Let $E$ be a Banach space, and $\Lambda$ an open set in $E$. Assign a subspace $M(x)$ in $E$ to each $x$ in $\Lambda$, but
the dimension of $M(x)$ can be infinite. In this
section, we consider the family $\mathcal{F}$ consisting of all $M(x)$
over $\Lambda,$ and investigate the sufficient and necessary
conditions for $\cal{F}$ being $c^1$ integrable at a point $x_0$ in
$\Lambda$.

{\bf Definition 2.1.}\quad Suppose that $M(x_0)$ is split in $E$, and $E=M(x_0)\oplus E_*.$  The set
$$J(x_0,E_*)=\{x\in\Lambda:M(x)\oplus E_*=E\}$$
is called the co-tailed set of $\mathcal{F}$ at $x_0$.

{\bf Theorem 2.1.}\quad{\it If $\mathcal{F}$ at $x_0$ is $c^1$ integrable, and let $S$ be its integral submanifold in $E$
, then

$(i)$ $\,M(x_0)$ is splitting in $E$;

$(ii) \,$ for an arbitrary complement $E_*$ of $ M (x_0),$ there exists a neighborhood $U_0$ at $x_0$, and a $c^1$ diffeomorphism $\varphi $ from $U_0$ onto $\varphi(U_0)$, such that
 $$\varphi'(x_0) =I_E, \, \, \varphi (S\cap U_0) \, \, \mbox{is an open set in $M(x_0),$}  $$ $$\varphi'(x) M(x)=M(x_0) \, \, \forall x\in S\cap U_0 ,\eqno(2.1)$$ and
 $$J(x_0, E_*)\supset U_0\cap S . \quad \eqno (2.2 )$$}

{\bf Proof:} \quad According to the definition of submanifolds in the Banach space $E,$ there
exists a subspace $E_0$ splitting in $E,$ a neighborhood $U_0$ at $x_0$, and a $c^1$ diffeomorphism $\varphi_1 $ from $U_0$ onto $\varphi_1(U_0)$, such that $\varphi_1 (S\cap U_0)$ is an open set in $E_0$. A key to proving the theorem is to verify the following formula:$$\varphi_1'(x) T_x S= E_0 \, \, \mbox{for} \, \, x\in S\cap U_0 . \quad \eqno (2.3)$$
Let $\delta $ be a sufficiently small positive number, $x$ a point in $S\cap U_0,$  and $c_{\delta,x } (t)$ from $(-\delta, \delta )$ to $S$ represent a $c ^ 1$ map that satisfies $c_{\delta,x } (0) =x.$ Since $\varphi_1 (S\cap U_0)$ is an open set in $E_0,$  $$\frac{d}{dt}(\varphi_1 \circ c_{\delta, x})\mid_{t=0}=\varphi_1'(x)\dot{c_{\delta, x}}(0)\in E_0 \, \, \mbox{for any $c^1$ curve} \, \, c_{\delta, x}(t).$$ By the definition of tangent space of $S$ at $x,$  $ T_x S=\{\dot{c_{\delta, x}}(0): \, \, \forall c^1 \, \, \mbox{curve} \, \, c_{\delta,x } (t)\},$
 $$\varphi_1'(x)T_x S\subset E_0.$$
 Conversely, since $\varphi_1 (S\cap U_0)$ is an open set in $E_0$,
for any $e\in E_0$ there exists a sufficiently small positive number $\delta$, such that for any $t\in (-\delta, \delta )$, $\varphi_1 (x)+t e \in \varphi_1 (S\cap U_0)$. Let
$ c_{\delta, x}(t)=\varphi_1 ^{-1}(\varphi_1 (x)+te ).$ Obviously, $ c_{\delta, x}(t)$ is
a $c^1$ curve in $S\cap U_0$, satisfying $c_{\delta, x}(0)=x.$ Directly,
$$\dot{c_{\delta, x}}(0)=(\varphi_1 ^{-1})'(\varphi_1(x))e=\varphi'(x)_1^{-1}e \, \, \forall e\in E_0,$$
and $$\varphi_1'(x)\dot{c_{\delta, x}}(0)=e \, \, \mbox{for any}  \,  e\in E_0.$$
Combined with the above results, we can confirm formula $(2.3).$ Note that since  $S \cap U _ 0$ is the integral submanifold of $\mathcal{F}$ at $x_0,$ $M(x)=T_x S \, \, \mbox{for any} \, \, x\in S \cap U _ 0.$ The formula $(2.3)$ derives  $$M(x) =\varphi_1'(x)^{-1}E_0 \,\, \mbox {for each} \, \, x \, \, \mbox{in} \, \, S\cap U_0 $$ In particular,
$$M(x_0)=\varphi_1'(x_0)^{-1}E_0, \quad \eqno (2.0)$$
We are now beginning to prove the theorem. Condition $(i)$ is valid because of the following two results : $E_0$ is split in $E$ and formula $(2.0).$
 Consider $$\varphi (x)=\varphi_1'(x_0)^{-1} \varphi_1(x) \, \, \forall x \in U_0.$$
Obviously, $\varphi'(x_0)=I_E$ and $\varphi'(x)=\varphi_1'(x_0)^{-1}\varphi_1'(x)$ for all x in $U_0.$
By the formula $(2.0)$ $\varphi_1'(x_0)\mid_{M(x_0)}$ belongs to $B^\times (M(x_0), E_0),$ where $B^\times (M(x_0), E_0)$ represent the set of all invertible operators in $B (M(x_0), E_0).$ So, according to the following three results: $ \varphi (U_0\cap S)= \varphi_1'(x_0)^{-1}(\varphi_1 (U_0\cap S)),$  $\varphi_1 (U_0\cap S)$ is an open set in $E_0,$ and the formula $(2.0)$, we infer that $\varphi (U_0\cap S)$ is an open set in $M(x_0)$, and for $x\in S\cap U_0,$
$$\varphi'(x) M(x)= \varphi_1'(x_0)^{-1} \varphi_1'(x) M(x)= \varphi'(x_0)^{-1}E_0=M(x_0).$$ This shows that formula $(2.1)$ is true.

The remaining proof of the theorem is to verify formula $(2.2).$
Let $E_*$ be an arbitrary complement of $ M (x_0).$  Consider the projection as follows,
   $$ P_x=\varphi'(x)^{-1}P^{E_*}_{M(x_0)}\varphi'(x)  \quad\quad\forall x\in S\cap U_0.$$
Clearly, $P^2_x=P_x$, $\lim\limits_{x\rightarrow x_0}P_x=P^{E_*}_{M(x_0)},$ and  $P_x$ is a generalized inverse of itself. By $(2.1)$
$$R(P_x)= \varphi'(x)^{-1} M(x_0)=M(x) .$$
We can now apply the operator rank theorem (Theorem $1.4$) to complete the proof of the theorem $2.1$. Taking $X$ and $T_x$ in Theorem $1.4$ as $S\cap U_0$ and $P_x$, respectively, the corresponding theorem shows that
there exists a neighborhood $S \cap U _ 1 $ at $x_0$ in $S \cap U _ 0,$ such that for $x\in S\cap U_1$, $$\|P_x-P_{x_0}\|< \|P_{M(x_0)}^{E_*}\|^{-1}, \, \, \mbox{and} \, \, R(P_x)\cap N(P_{x_0})=\{0\},$$ where $ U _ 1$ is a neighborhood at $x_0$ in $U_0.$
 By the equivalence of $(i)$ and $(iii)$ in Theorem $1.1$,
$$R(P_x)\oplus N(P_{M(x_0)}^{E_*})=E, \, \, \mbox{i.e.,} \, \, M(x)\oplus E_*=E \, \, \mbox{for} \, \, x\in S\cap U_1. $$
 This means
 $J(x_0, E_*)\supset U_0\cap S.$  Without lose of generality, still write $U_1 $ as $U_0.$ Formula $(2.2)$ is proved. The proof ends \quad $\Box$

Theorem $1.5$ shows that the $M (x)$ for $x \in J (x_0, E_ *)$ has a unique operator valued coordinate $\alpha (x)$ that satisfies $$M(x)=\{e+\alpha (x)e:\forall e\in M (x_0)\}.$$ Let $C^1(V, E_*)$ denote all of $c^1$ maps from an open set $V$ in $M(x_0)$ to subspace $E_*,$ and $C_0^1(V, E_*)=\{\psi \in C^1(V, E_*): \psi (P_{M(x_0)}^{E_*}x_0)=P_{E_*}^{M(x_0)}x_0 \}.$

{\bf Theorem 2.2}\ (Frobenius theorem).\quad{\it

If $\mathcal{F}=\{M(x):x\in \Lambda \}$ is $c^1$
integrable at $x_0$, then

(i) $M(x_0)$ is splitting in $E$, and

(ii) for each complement $E_*$ of $ M (x_0),$
there exists a neighborhood $V$ at $P_{M(x_0)}^{E_*}x_0$ in $M(x_0)$, and a map $\psi$ in $C^1(V, E_*),$  satisfying both that
$$ J (x_0, E_ *) \supset \{x+\psi (x): \forall x\in V\},\quad \eqno (2.4)$$ and for each $x\in V$, $$\psi'(x)=\alpha(x+\psi(x)), \, \, \mbox{and} \, \,  \psi(P_{M(x_0)}^{E_*}x_0)=P^{M(x_0)}_{E_*}x_0 .\quad\eqno(2.5) $$
 Where  $\psi'(x)$ is the Fr\'{e}chet derivative of  $\psi$ $at$ $x.$ \\
( The decomposition $E=M(x_0)\oplus E_*$  brings a coordinate $(x,\psi (x))$ of $x+\psi (x),$ so, $\alpha(x+\psi(x))$is also can write as $\alpha(x,\psi(x))$ .)

Conversely, if the condition $(i)$ holds, and there is a complement $E_*$ of $M(x_0)$, a neighborhood $V$ at $P_{M_0}^{E_*}x_0$ in $M(x_0)$, and a map $\psi \in C^1(V, E_*),$ satisfying the equalities $(2.4)$ and $(2.5)$, then $\mathcal{F}$ is $c^1$ integrable at $x_0.$ }

{\bf Proof:}\quad Suppose that $\mathcal{F}$ is $c^1$-integrable at $x_0$ and $S$ is its integral submanifold at $x_0.$
We want to verify conditions $(i)$ and $(ii)$ in the theorem.
According to Theorem $2.1$, the conditions $(i)$ in Theorem $2.2$ is true, and there is the neighborhood $U_0$ at $x_0,$ and  the $c^ 1$ diffeomorphism $\varphi$ with $\varphi'(x_0)=I_E $, such that the set $\varphi (S \cap U_0)$ is an open set in $M (x_0),$ containing the point $\varphi (x _0).$ For simplicity, write $V_0=\varphi (S\cap U_0)$ in the sequel. Next, we use $\varphi$ and $V_0$ to determine the neighborhood $V$ and the map $\psi $ in the condition $(ii).$
For any complement $E_*$ of $ M (x_0),$ consider
$$ \varphi _0=P^{E_*}_{M(x_0)}\varphi ^{-1}\mid_{V_0} : V_0\rightarrow M(x_0),$$ and $$\varphi _1=P_{E_*}^{M(x_0)}\varphi ^{-1}\mid_{V_0} : V_0\rightarrow E_*.$$
obviously, $$ \varphi_0 (\varphi (x_0))=P_{M(x_0)}^{E_*}x_0 ,$$ and since $ \varphi'(x_0)=I_E,$
$$\varphi _0'(\varphi(x_0))=P^{E_*}_{M(x_0)}(\varphi ^{-1})'(\varphi (x_0))=P^{E_*}_{M(x_0)}\varphi'(x_0)^{-1}=I_{M(x_0)} .$$ Thus, by the inverse mapping theorem, there exists in $V_0$ a neighborhood at $\varphi (x_0)$, without lose of generality , still write it as $V_0$, such that
 $ \varphi _0 $ from $ V_0$ onto $\varphi _0(V_0)$ is a $c^1$ diffeomorphism. Let $ V=\varphi _0 (V_0), \, \, \mbox{and} \, \, \psi=\varphi _1\circ  \varphi _0 ^{-1}: V\rightarrow E_*.$ We claim that they meet the condition $(ii).$ Clearly,
 \begin{eqnarray*} S\cap U_0 &=&\varphi ^{-1}(V_0)\\ &=& \varphi ^{-1}(\varphi ^{-1}_0(V))\\&=& (P^{E_*}_{M(x_0)}\varphi^{-1}+P_{E_*}^{M(x_0)}\varphi ^{-1} ) (\varphi ^{-1}_0(V))\\ &=& (\varphi _0+P_{E_*}^{M(x_0)}\varphi ^{-1})(\varphi ^{-1}_0(V))\\&=& (I_{M(x_0)}+ \psi)(V).\end{eqnarray*}
Hence,
  $$S\cap U_0=\{x+\psi (x): \forall x\in V\}, \quad \eqno (2.6)$$
 and by $(2.2)$, $\psi $ meets $(2,4)$.Next, we want to verify that $\psi $ meets $(2.5)$. Since the above equality $ \varphi_0 (\varphi (x_0))=P_{M(x_0)}^{E_*}x_0 ,$ therefore
\begin{eqnarray*} \psi(P_{M (x_0)}^{E_*}x_0) &=& (\varphi _1\circ \varphi _0 ^{-1})(P_{M(x_0)}^{E_*}x_0)=\varphi _1(\varphi _0 ^{-1}(P_{M_0}^{E_*}x_0))\\ &=&\varphi _1(\varphi(x_0))= P_{E_*}^{M(x_0)}\varphi ^{-1}(\varphi(x_0))=P_{E_*}^{M(x_0)}x_0 .\end{eqnarray*}
In the other hand, the condition $(2.4)$ shows that for $x\in V, \, x+\psi (x)\in J (x_0, E_ *),$ so that the theorem $1.5$ can be applied to $M(x+\psi (x))$ and $M(x_0).$ Therefore, $$M(x+\psi(x))=\{e+\alpha(x+\psi (x))e : \forall e\in M_0\} \, \, \mbox{for} \, \, x\in V. \quad \eqno (2.7) $$
 Since $S\cap U_0$ is the integral submanifold of $\mathcal{F}$ at $x_0,$ $$T_{x+\psi (x)}(S\cap U_0)=M(x+\psi(x)) \, \, \forall x\in V .$$ Therefore, $$T_{x+\psi (x)}(S\cap U_0)=\{e+\alpha(x+\psi (x))e : \forall e\in M_0\}.$$
 To complete the proof of the necessity of the theorem we need to calculate $T_{x+\psi (x)}(S\cap U_0).$

 Let $c_{\delta, x+\psi (x)}(t)$ from $(-\delta, \delta)$ to $S\cap U_0$ represent $c^1$ curve via the point $(0, x+ \psi(x)).$ By $(2.6)$ it can be written as $z(t)+\psi (z(t))$, where $z(t)$ from $(-\delta, \delta)$ to  $M(x_0)$ is $c^1$ curve through the point $(0, x)$.  We have that for any $c^1$ curve $c_{\delta, x+\psi (x)}(t)$ from $(-\delta, \delta)$ to $S\cap U_0$,  $$\dot {c}_{\delta, x+\psi (x)}(0)=\dot{z}(0)+\psi'(x)\dot{z}(0).$$ therefore, for all $x\in V$
 $$T_{x+\psi (x)} (S\cap U_0)=\{e +\psi'(x)e : \forall e\in M(x_0)\}. \quad \eqno (2.8)$$
Combining $(2.7)$ and $(2.8)$ we deduce
$$\psi'(x)=\alpha (x+\psi (x)) \, \, \mbox{for} \, \, x\in V.$$
The necessity of the theorem is proved.

Conversely, suppose that $M(x_0)$ splits $E,$ and for a complement $E_*$ of $M(x_0),$ there is a neighborhood $V$ at $P_{M(x_0)}^{E_*}x_0$ in $M(x_0)$ and a map $\psi \in C^1(V, E _ *)$ satisfying that $(2.4)$ and $(2.5).$ We want to show that $\mathcal {F}$ at $x_0$ is $c^1$ integrable. Consider
$$V^*=\{x\in E:P^{E_*}_{M(x_0)}x\in V\} ,$$ and $$ \Phi(x)=x+\psi(P^{E_*}_{M(x_0)}x)\, \, \mbox{for } \, \, x\in V^*.$$
Let $$S=\{x+\psi (x): \forall x\in V\}.$$ Clearly, $V^*$ is an open set containing $V$ and $S$, and  $$\Phi^{-1}(z)=z-\psi(P^{E_*}_{M(x_0)}z) \, \, \mbox{for} \, \, z\in V^*.$$
Indeed,
$$(\Phi^{-1}\circ \Phi )(z)=\Phi (z)-\psi(P_{M(x_0)}^{E_*}\Phi(z))=z+\psi(P_{M(x_0)}^{E_*}z)- \psi(P_{M(x_0)}^{E_*}z)=z \, \, \forall z\in V^*,$$
 and
$$(\Phi \circ \Phi^{-1})(z)=Y(z)+\psi (P_{M(x_0)}^{E_*}Y(z))=z-\psi(P_{M(x_0)}^{E_*}z)+\psi(P_{M(x_0)}^{E_*}z)=z \, \, \forall z\in V^*.$$
Therefore, $\Phi $ from $V^*$ onto itself is $c^1$ diffeomorphism. Now, for each $x\in S,$ we have the triples $(\Phi^{-1}, V^*, E )$ and $( \Phi ^{-1}\mid_S, V^*\cap S, M(x_0 )).$ Obviously, $\Phi^{-1}(V^* \cap S)=\Phi^{-1}(S)=V ,$ so it is an open set in $M(x_0).$ Hence
 $S$ is the $c^1$ integral submanifold with an atlas consisting of single coordinate chart $( \Phi ^{-1}\mid_S, V^*\cap S, M(x_0 )).$
 To complete the proof we now only need to verify that $S$ is tangent to $M(y)$ at any point $y\in S,$ i.e., $T_y S= M(y).$

 Let $c_{\delta,y} (t)$ from $(-\delta,\delta)$ to $S$ be any  $c^1$ curve through the point $(0, y)$, where $\delta $ is a small enough positive number. According to the definition of $T_yS$ this requires a proof of
 $$M(y)=\{\dot{c_{\delta,y}}(0):\forall \, c^1 \, \, \mbox{curve} \, \, c_{\delta,y}(t)\}.$$
Because $S=\{x+\psi (x): \forall x\in V\}$ we can write $$c_{\delta, y} (t)= P_{M(x_0)}^{E_*} c_{\delta, y}(t) + \psi (P_{M(x_0)}^{E_*} c_{\delta, y} (t)) \in S \, \, \mbox{for} \, \, t\in (-\delta, \delta),$$
and derive $$\dot{c_{\delta, y}}(0)= P_{M(x_0)}^{E_*}\dot {c_{\delta,y}}(0)+\psi'(P_{M(x_0)}^{E_*}y)P_{M(x_0)}^{E_*}\dot{c_{\delta x,y}}(0).$$
 By $(2.5)$
\begin{eqnarray*} \dot{c_{\delta, y}}(0)&=& P_{M(x_0)}^{E_*}\dot{c_{\delta, y}}(0)+\alpha (P_{M(x_0)}^{E_*}y+ \psi(P_{M(x_0)}^{E_*}y))P_{M(x_0)}^{E_*}\dot{c_{\delta,y}}(0)\\ &=&P_{M(x_0)}^{E_*}\dot{c_{\delta, y}}(0)+\alpha (y)P_{M(x_0)}^{E_*}\dot{c_{\delta, y}}(0),
\end{eqnarray*}
 this indicates that $\dot{c_{\delta, y}}(0)\in M(y),$ and hence, $T_y S \subset M(y).$ Conversely, for each $e\in M(y)$ let $e=e_0+\alpha (y) e_0 $ for $e_0\in M(x_0)$, and $y=x+\psi (x)$ for $x\in V.$ Consider
 $$ c_{\delta, y}(t)=( x+t e_0 )+ \psi (x+t e_0 ) \, \, \mbox{for} \, \, t\in (-\delta, \delta).$$
Clearly, $c_{\delta, y}(t)\in S$ for $t\in (-\delta, \delta)$ and $c_{\delta, y}(0)=y.$ By $(2.5)$
\begin{eqnarray*} e &= & e_0+\alpha (y) e_0\\&=& e_0+\alpha (x +\psi (x))e_0 \\&=&e_0+ \psi'(x) e_0 \\ &=& \dot{c_{\delta, y}}(0)\in T_y S \end {eqnarray*}
 This shows $T_y S\supset M(y).$ Combining the two results above, we conclude that $M(y)=T_y S \, \, \forall y \in S.$ The proof ends.\quad$\Box$

\vskip 0.2cm\begin{center}
{\bf 3\quad Frobenius theorem  with a Trivial Co-Tailed Set }\end{center}\vskip 0.2cm

Let $E$ be a Banach space, $U$ an open set in $E$, and  $\mathcal{F} \, \,=\{M(x)\}_{x\in U}$ be a family of subspaces in $E$.

\textbf{Definition $3.1$.}\quad Assume that $x_0 \in U_0$ and $M(x_0)\oplus E_*=E.$ The co-tailed set $J(x_0,E_*)$ is said
to be trivial provided  $x_0$ is an inner point of $J(x_0,E_*).$

 \textbf { Theorem $3.1$.} \quad {\it If $M(x_0)\oplus E_*=E$ and $J(x_0,E_*)$ is trivial, then the $c^1$ integrability of $\mathcal{F}$ at $x_0$ and the solvability of the differential equation $(2.5)$ with the initial value $\psi (P_{M(x_0)}^{E_*}x_0)=P^{M(x_0)}_{E_*}x_0$ are equivalent.}

$\quad ${\bf Proof:}\quad Assume that $J(x_0,E_*)$ is trivial, and let $W_0$ be a neighborhood of $x_0$ in $J(x_0, E_*).$
Consider the mapping $$x+y \, \, \mbox{from $M(x_0)\times E_*$ into $E$ .}$$ Since the mapping $x+y$ is continuous at point $ (P_{M(x_0)}^{E_*}x_0, P_{E_*}^{M(x_0)}x_0)$, there exists a neighborhood $U_0$ of $P_{M(x_0)}^{E_*}x_0$ in $M(x_0)$, and a neighborhood $V_0$ of $P_{E_*}^{M(x_0)}x_0$ in $E_*$, such that $$x+y \in W_0 \, \, \mbox{for all} \, \,(x, y)\in U_0 \times V_0.$$
Obviously, for each $\psi \in C^1_0(U_0, E_*),$
there is a neighborhood $V$ at $P_{M(x_0)}^{E_*}x_0$ in $U_0$, such that $\psi (V)\subset V_0.$ This causes $$\{ x +\psi (x): \, \forall x\in V\} \subset J(x_0,E_*).$$
Therefore, when $J(x_0, E_*)$ is trivial, for any $\psi \in C_0^1(U_0, E_*)$ there is  neighborhood $V$ at the point $P_{M(x_0)}^{E_*}x_0$ in $U_0$ ($V$ depends on $\psi$ ), making that $\{x+\psi (x): \, \forall \in V \} \subset J(x_0, E_*).$ Now, we conclude that in the case of trivial $J(x_0, E_*)$, the two conditions $(i)$ and $(ii)$ in Theorem $2.2$ reduce to the solvability of the differential equation $(2.5)$ with the initial value
$\psi (P_{M(x_0)}^{E_*}x_0)=P^{M(x_0)}_{E_*}x_0.$ The proof ends.\quad $\Box$

$ \textbf {Corollary $3.1$} :$ suppose that $E$ is an Euclidean space, and $J(x_0,E_*)$ is trivial. Then the $c^1$ integrability of $\mathcal{F}$ at $x_0$ and the solvability of the differential equation $(2.5)$ with the initial value $\psi (P_{M(x_0)}^{E_*}x_0)=P^{M(x_0)}_{E_*}x_0$ are equivalent.

The following example illustrates Theorems $3.1$ and Theorem $2.2$ and their use, despite its simplicity.

 {\bf Example}.\quad Let $E=\mathbf{R}^3,
  \Lambda=\mathbf{R}^3\setminus(0,0,0)$ and
  $$M(x_1,y_1,z_1)=\{(x,y,z)\in\mathbf{R}^3 : x_1x+y_1y+z_1z=0\}\, \,\mbox{for all} \, \,
  (x_1,y_1,z_1)\in\Lambda.$$
  Consider the family of subspaces $\mathcal{F}=\{M(x_1,y_1, z_1):\forall (x_1,y_1, z_1)\in \Lambda\}$. Apply Corollary $3.1$ to
  determine the integral surface of $\mathcal{F}$ at $(0,0,1)$.

Let $U_0=\{(x, y, z)\in\mathbf{R}^3: z >0\},$ and let $E_*=\{(0,0, z)\in \mathbf{R}^3:\forall z\in R\}$.
 It is easy to see $ J ((0,0,1), E_ *)\supset U_0.$ In fact, $ M (x_1, y_1, z_1) \cap E_ * = \{(0,0,0)\}$ for any $(x_1, y_1, z_1)\in U_0,$  and so, $ M(x_1,y_1,z_1) \oplus E_*= \mathbf{R}^3$ because $dim M(x_1,y_1,z_1)=2.$ This produces $J ((0,0,1), E_ *) \supset U_0.$ Therefore, point $(0,0,1)$ is an inner point of  $J((0,0,1),E_*)$, that is, $J( (0,0, 1), E_*)$ is trivial. Directly,
 \begin{eqnarray*}
  M(x_1,y_1,z_1)&=&\{(x,y,-\frac{x x_1+yy_1}{z_1})): \forall (x,y)\in R^2\}\\
 & =&\{(x,y,0)+(0,0,-\frac{xx_1+yy_1}{z_1}): \, \forall (x,y)\in R^2\}\\ &=& \{(x,y,0)+\alpha (x_1,y_1,z_1)(x,y,0): \, \forall (x,y)\in R^2\},\end{eqnarray*}
where $$\alpha (x_1,y_1,z_1)(x,y,0)= (0,0,-\frac{xx_1+yy_1}{z_1})= (0,0, - \frac{1}{z_1}(x_1,y_1)\bullet(x, y)),$$ and $(x_1,y_1)\bullet(x, y)=xx_1+yy_1.$
Note $M(0,0,1)=\{(x, y, 0): \forall (x, y) \in R^2\},$ and write it as $M_0.$
Obviously,  $P_{M_0}^{E_*}(0,0,1)=(0,0,0)$ and $P^{M_0}_{E_*}(0,0,1)=(0,0,1).$ Now let's solve the differential equation $(2.5).$ Consider a map $\psi \in C_0^1 ( M_0, E_*).$ Since $\psi'(x_1,y_1,0)$ belongs to $B (M_0,E_*),$
\begin{eqnarray*}
\psi'(x_1,y_1,0)(\Delta x_1,\Delta y_1)&=&(\frac{\partial\psi}{\partial x_1}(x_1,y_1,0),\frac{\partial\psi}{\partial y_1}(x_1,y_1,0))\bullet(\Delta x_1,\Delta y_1)\\ &=&-\frac{(x_1,y_1)}{\psi(x_1,y_1,0)}\bullet (\Delta x_1,\Delta y_1).\end{eqnarray*} Thus
 $${\psi(x_1,y_1,0) \frac{\partial\psi}{\partial x_1}(x_1,y_1,0)=-x_1} \, \, \mbox{and} \, \, {\psi (x_1,y_1,0) \frac{\partial\psi}{\partial y_1}(x_1,y_1,0)=-y_1},$$ so,
$$\frac{\partial\psi ^2}{\partial x_1}(x_1,y_1,0)=-\frac{1}{2}x_1 \, \mbox{and} \, \, \frac{\partial\psi ^2}{\partial y_1}(x_1,y_1,0)=- \frac{1}{2}y_1 .$$
Through integration,
$$\psi ^2 (x_1,y_1,0)=-x_1^2 + c(y_1) \, \, \mbox{for a } \, \, c(y_1)\in C^1(-1,1),$$ and then
$c'(y_1)=-\frac{1}{2}y_1 .$ Therefore, $c(y_1)=-y_1^2 +c$ where $c$ is a constant. Additionally, since $ \psi(0,0,0)=(0,0,1)$ we have $\psi^2(x_1,y_1,0)=1-x_1^2-y_1^2$. We can now say that $\{(x,y, \sqrt{1-x^2-y^2}): \mbox {for both} \, \|x\|,\|y\|\leq 1\}$ is the integrated surface of $\mathcal{F}$ at $(0,0,1).$

\textbf{Theorem $3.2$.}\quad{\it Suppose that the set $U$ in $E$ is an open set containing the point $x_0$, and $f$ is a $c^1$ mapping from  $U$ to Banach space $F.$ If $x_0$ is a generalized regular point of $f$ then $\mathcal{F}=\{N(f'(x)):\forall x\in U_0\}$ is $c ^ 1$ integrable at $x_0,$ and the corresponding $J (x_0, E_ *)$ is trivial. }

\textbf{Proof:}\quad Since $x_0$ is a generalized regular point of $f$, Definition $1.2$ and Theorem $1.3$ show that for any generalized inverse
$f'(x_0)^+$ of $f'(x_0)$, there exists a neighborhood $U_0$  at $x_0$ such that
$$R(f'(x))\cap N(f'(x_0)^+)=\{0\}\quad\mbox{for all} \, \, x\in U_0.$$
Since $f'(x)$ is continuous at $x_0$, we can assume that for $x\in U_0$,
$$\|f'(x)-f'(x_0)\|<\|f'(x_0)^+\|^{-1} \, \, \, \mbox{and} \, \, \, R(f'(x))\cap N(f'(x_0)^+)=\{0\}.\eqno(3.1)$$
Let $E_*=R(f'(x_0)^+),$ and write $N(f'(x))$ and $N(f'(x_0))$ as $N_x$ and $N_0$, respectively. We first go to show that $x_0$ is an inner point of
$J(x_0,E_*)$. Since $x_0$ is a generalized regular point of $f$, the two conclusions in $(3.1)$ and the equivalence of conditions $(i)$ and $(iv)$ in Theorem $1.1$ guarantee  a neighborhood $U_0$ at $x_0,$ satisfying that for any  $x \in U_0,$ $$N_x\oplus E_*=E \, \, \, \mbox{and} \, \,\, \|f'(x)-f'(x_0)\|<\|f'(x_0)^+\|^{-1} .$$
This shows that $J(x_0,E_*)\supset U_0,$ i.e., a trivial $J (x_0, E_ *)$.
Next, we will determine $\alpha (x)$ about $x\in U_0.$ Obviously, for any $e\in N_x,$
$$P^{E_*}_{N_x}P^{E_*}_{N_0}e=P^{E_*}_{N_x}(P^{E_*}_{N_0}e+P^{N_0}_{E_*}e)=P^{E_*}_{N_x}e=e ,$$ and for any  $x\in N_0,$
$$P^{E_*}_{N_0}P^{E_*}_{N_x}x=P^{E_*}_{N_0}(P^{E_*}_{N_x}x+P^{N_x}_{E_*}x)=P^{E_*}_{N_0}x=x.$$
 So for any  $e\in N_x,$
$$e=P^{E_*}_{N_0}e+P^{N_0}_{E_*}e=P^{E_*}_{N_0}e+P^{N_0}_{E_*}P^{E_*}_{N_x}P^{E_*}_{N_0}e .$$
We can now speculate on the result as follows,$$ \left.\alpha(x)=P^{N_0}_{E_*}P^{E_*}_{N_x}\right|_{N_0} \, \,  \mbox{for all} \, \, x\in U_0.$$
For the proof of it we only need to verify that
$$ P^{E_*}_{N_0} N_x= N_0 \, \,  \forall x \in U_0.$$
 Because $J(x_0,E_*)$ contains $U_0,$  $$ P^{E_*}_{N_x}e_0=e_0 - P_{E_*}^{N_x}e_0 \, \, \mbox{for any $x\in U_0$} \, \, \mbox{ and any $e_0\in N_0$,}$$
and let $e=P^{E_*}_{N_x}e_0 ,$ then $e$ belongs to $N_x$ and satisfies $$P^{E_*}_{N_0}e=P^{E_*}_{N_0}e_0=e_0.$$ This shows $ P^{E_*}_{N_0} N_x=N_0 \, \,  \forall x \in U_0.$
 According to theorem $1.5,$ the operator valued coordinate $\alpha (x)$ of $N_x$ is unique, hence we assert that
$$ \left.\alpha(x)=P^{N_0}_{E_*}P^{E_*}_{N_x}\right|_{N_0} \, \,  \mbox{for all} \, \, x\in U_0.$$
Let $f'(x_0)$ and $f'(x)$ replace $A$ and $T$ in theorem $1.1,$ respectively. Then, the two conclusions in $(3.1)$ guarantee the establishment of the corresponding theorem for $x\in U_0.$  Therefore, by the equivalence of the conditions $(i)$ and $(ii)$ in Theorem $1.1$, for any $ x\in U_0,$ $f'(x)$ has the following generalized inverse $f'(x)^+$:
$$ f'(x)^+= f'(x_0)^+ C^{-1}_{f'(x_0)}(f'(x_0)^+,f'(x))=D^{-1}_{f'(x_0)}(f'(x_0)^+, f'(x))f'(x_0)^+ ,$$ satisfying that for $x\in U_0,$ $$N(f'(x)^+)=N(f'(x_0)^+), \, \, \mbox{and} \, \, R(f'(x)^+)=R(f'(x_0)^+)= E_*.$$
 Therefore, for $x\in U_0,$

 $$P^{E_*}_{N_x}= P^{R(f'(x)^+)}_{N_x} = I_E-f'(x)^+f'(x) = I_E-D^{-1}_{f'(x_0)}(f'(x_0)^+,f'(x))f'(x_0)^+f'(x).$$
 Thus,
\begin{eqnarray*}\alpha(x)&=&\left.P^{N_0}_{E_*}P^{E_*}_{N_x}\mid_{N_0}\right.\\
&=&P^{N_0}_{E_*}\left.(I_E-D^{-1}_{f'(x_0)}(f'(x_0)^+, f'(x))f'(x_0)^+f'(x))\right.\mid_{N_0}\\
&=&P^{N_0}_{E_*}D^{-1}_{f'(x_0)}(f'(x_0)^+,f'(x))\left.(D_{f'(x_0)}(f'(x_0)^+,f'(x))-f'(x_0)^+f'(x)\right.\mid_{N_0}\\
&=&P^{N_0}_{E_*}D^{-1}_{f'(x_0)}(f'(x_0)^+,f'(x))\left.(I_E-f'(x_0)^+ f'(x_0)\right.\mid_{N_0}\\
&=&P^{N_0}_{E_*}D^{-1}_{f'(x_0)}(f'(x_0)^+,f'(x))P^{E_*}_{N_0} \quad \mbox{for all} \, \, x\in U_0,
\end {eqnarray*} this indicates that
$$\alpha(x)=P^{N_0}_{E_*}D^{-1}_{f'(x_0)}(f'(x_0)^+,f'(x))P^{E_*}_{N_0} \quad \mbox{for all} \, \, x\in U_0. \eqno(3.2)$$
Evidently,
\begin{eqnarray*}
  (f'(x_0)^+(f(x)-f(x_0))+P^{E_*}_{N_0}x)'&=&f'(x_0)^+f'(x)+P^{E_*}_{N_0}\\&=&
f'(x_0)^+(f'(x)-f'(x_0))+f'(x_0)^+f'(x_0)+P^{E_*}_{N_0}\\&=&I_E+f'(x_0)^+(f'(x)-f'(x_0))\\&=&D_{f'(x_0)}(f'(x_0)^+,f'(x))
.\end {eqnarray*}  Let $\varphi(x)=f'(x_0)^+(f(x)-f(x_0))+P^{E_*}_{N_0}x,$ then
\begin{eqnarray*}
\varphi'(x)&=&f'(x_0)^+f'(x)+P_{E_*}^{N_0}=f'(x_0)^+f'(x)+(I_E-f'(x_0)^+f'(x_0))\\ &=&I_E+f'(x_0)^+(f'(x)-f'(x_0))=D_{f'(x_0)}(f'(x_0)^+,f'(x)).
\end {eqnarray*}So
$$\varphi'(x)=D_{f'(x_0)}(f'(x_0)^+,f'(x)) \, \, \mbox{for all} \, \, x\in U_0 .\eqno(3.3)$$
Obviously, $\varphi (x_0)= P_{N_0}^{E_*}x_0,$ and $\varphi'(x_0)=I_E.$
  Thus,
by the inverse mapping theorem, there exists a neighborhood at $x_0$, without loss of generality, still write it as $U_0$, such that $\varphi$ from $U_0$ onto $\varphi (U_0)$ is a $c^1$ diffeomorphsm, where $\varphi (U_0)$ is an open set containing $P_{N_0}^{E_*}x_0,$ and write it as $V_0$. Then, by the equalities $(3.2)$ and $(3.3),$
\begin{eqnarray*}
\alpha(x)&=&P^{N_0}_{E_*}D^{-1}_{f'(x_0)}(f'(x_0)^+,f'(x))|_{N_0}\\
&=&P_{E_*}^{N_0}\varphi'(x)^{-1}|_{N_0}\\
&=&P_{E_*}^{N_0}(\varphi^{-1})'(y)|_{N_0}. \end{eqnarray*}
This points to the following key equality:
$$\alpha(\varphi^{-1}(y))=P_{E_*}^{N_0}(\varphi^{-1})'(y)|_{N_0} \, \, \mbox{for all} \, \, y\in V_0, \eqno(3.4)$$
which reveals how to derive the solution $\psi$ of the differential equation $(2.5).$ Now, we are going to derive $\psi$ satisfying $ \varphi^{-1}(z)=z+\psi(z).$
 Consider the $c^1$ map as follows,
$$\varphi _0=P_{N_0}^{E_*}\varphi ^{-1} : V_0 \cap N_0\rightarrow N_0 . $$
From $\varphi(x_0)=P_{N_0}^{E_*} x_0$ it is to see $\varphi_0(P_{N_0}^{E_*}x_0)=P_{N_0}^{E_*} x_0,$  that is, $P_{N_0}^{E_*} x_0$ is a fixed point of $\varphi_0$ , and
$$\varphi _0'(P_{N_0}^{E_*}x_0)=P_{N_0}^{E_*}(\varphi^{-1})'(P_{N_0}^{E_*}x_0)=P_{N_0}^{E_*}\varphi'(x_0)^{-1}=I_{N_0}.$$
Then by the inverse mapping theorem, there exists in $V_0\cap N_0$ a neighborhood at $P_{N_0}^{E_*}x_0,$ without loss of generality, still write it as $V_0\cap N_0,$ such that
$$\varphi _0 \, \mbox{from} \, V_0\cap N_0 \, \mbox{onto} \, \varphi _0(V_0\cap N_0) \, \, \mbox{is a $c^1$ diffeomorphism}.$$  Set $V=\varphi _0(V_0\cap N_0).$ It is easy to see that $V$ is an open set in $N_0$, containing the point $P_{N_0}^{E_*}x_0.$
We claim that
$$\varphi'_0(y)=I_{N_0} \, \, \forall y\in V_0\cap N_0.$$
 For any $y \in V_0\cap N_0,$ let $ x \in U_0$ satisfy $y =\varphi (x).$ Note $R(f'(x_0)^+)= E_*,$ then for any $e\in N_0,$
\begin{eqnarray*} \varphi'_0(y)e &=& P^{E_*}_{N_0}(\varphi^{-1})'(y)e= P^{E_*}_{N_0}\varphi'(x)^{-1}e\\&=& P^{E_*}_{N_0}D_{f'(x_0)}^{-1}(f'(x_0)^+,f'(x))e \, \, \mbox{by} \,\, (3.3)\\
&=& P_{N_0}^{E_*}(P_{N_0}^{E_*}+f'(x_0)^+f'(x) )D_{f'(x_0)}^{-1}(f'(x_0)^+,f'(x))e  \\ &=& P_{N_0}^{E_*}D_{f'(x_0)}(f'(x_0)^+,f'(x))D_{f'(x_0)}^{-1}(f'(x_0)^+,f'(x))e\\&=&e.\end {eqnarray*} This indicates that $\varphi'_0(y)=I_{N_0} \, \, \forall y\in V_0\cap N_0,$
 and hence, $$(\varphi_0^{-1})'(z)=I_{N_0} \, \  \forall z\in V. \quad \eqno (3.5)$$
Let $\varphi _1$ represent the map $P^{N_0}_{E_*}\varphi ^{-1}$ from $V$ to $E_*,$ and $\psi$ be the $c^1$ map as follows, $$\psi=\varphi _1\circ \varphi _0^{-1} \, \, \mbox{from $V$ to $E_*$}.$$
We now want to verify that $\psi $ is the solution of the initial valued problem $(2.5)$. Directly,
 \begin{eqnarray*}
\psi(P^{E_*}_{N_0}x_0)&=&\varphi_1\left(\varphi^{-1}_0(P^{E_*}_{N_0}x_0)\right)\\
&=&\varphi_1(P^{E_*}_{N_0}x_0)=P_{E_*}^{N_0}\varphi^{-1}(P^{E_*}_{N_0}x_0)\\
&=&P^{N_0}_{E_*}x_0 .
\end{eqnarray*}
In addition, let $y\in V_0\cap N_0$ satisfy $z=\varphi_0(y).$ Then for any $e\in N_0,$
\begin{eqnarray*}\psi'(z)e &=&P_{E_*}^{N_0}(\varphi^{-1}\circ \varphi_0^{-1})'(z)e\\&=& P_{E_*}^{N_)}(\varphi^{-1})'(\varphi_0^{-1}(z))\cdot(\varphi_0^{-1})'(z)e\\&=& P_{E_*}^{N_0}(\varphi^{-1})'(\varphi_0^{-1}(z))e \, \, \mbox{by (3.5)}\\&=& P_{E_*}^{N_0}(\varphi^{-1})'(y)e\\ &=&\alpha (\varphi^{-1}(y))e \, \, \mbox{by (3.4)},
\end {eqnarray*} and for each $e\in N_0,$$$\alpha(\varphi^{-1}(y))e= \alpha (\varphi_0(y)+ \varphi_1 (y))e=\alpha(z+\psi (z))e \, \, \mbox{for} \, \, z\in V. $$
Therefore, for all $z\in V,$ $$\psi'(z)=\alpha (z+\psi (z)).$$
Theorem $3.2$ is proved.\quad$\Box$

{\bf Note:} It is easy to see from Theorem $1.2$ that the following points are also generalized regular points: regular points, submersion points, immersion points, subimmersion points, Fredholm points,and semi-Fredholm points, so Theorem $3.2$ gives a large class of solvable initial valued problem of differential equations with geometric meaning. This is an interesting addition to geometrical methods in the theory of ordinary differential equations (see [An]). It's worth digging deeper.

\vskip 0.2cm\begin{center}{\bf 4\quad A Family of Subspaces with
Non-Trivial Co-Tailed Set and Smooth Integral
Submanifolds}\end{center}\vskip 0.2cm

Let $\mathbf{M}(X)=\{T\in
B(E,F):TN(X)\subset R(X)\}$ for $X\in B(E,F)$, which appears at the first time in [Caf]. In this section we consider the family of subspaces
 $\mathbf{\cal {F}}=\{\mathbf{M}(X)\}_{X\in B(E, F)}.$ First, let's introduce an example below.

{\bf Example}.\quad Consider the space $B(R^2)$  consisting of
all real $2\times 2$ matrices. Let $A$ and $\mathbf{E}_*$ be
$$\left(\begin{array}{cc}1&0\\ 0&0\end{array}\right)\quad{\rm
and}\quad\left\{\left(\begin{array}{cc}0&0\\
0&t\end{array}\right):\forall t\in R\right\},$$ respectively.Obviously
$$N(A)=\{(0,x):\forall x\in R\} \, ,  R(A)=\{(x,0),\forall x\in R\} , $$
$$M(A)=\left\{\left(\begin{array}{cc}t_{1 1}&t_{1 2}\\
t_{2 1}&0\end{array}\right),\quad\quad \forall t_{1 1}, \, t_{1 2} \, \, and \, \, t_{21} \in R  \right\}, $$
and $$M(A)\oplus\mathbf{E}_*=B(R^2).$$
Now we are going to verify that $J(A,\mathbf{E}_*)$ is nontrivial.Consider
$$A_\varepsilon=\left(\begin{array}{cc}1&0\\
0&\varepsilon\end{array}\right),\quad\quad\varepsilon\not=0.$$  Obviously $$ N(A_\varepsilon ) =\{0\} \, \, \mbox{and} \, \, M(A_\varepsilon ) =B(R^2),$$ therefore, $$\mbox{dim} M(A_\varepsilon)= 4, \mbox{and} \, \,
\lim\limits_{\varepsilon\rightarrow 0}A_\varepsilon=A.$$
Obviously, for any non zero $\varepsilon,$ $ \mbox{dim} M(A_{\varepsilon})=4,$ so $M(A_{\varepsilon})$ does not satisfy $M(A_{\varepsilon})\oplus \mathbf{E}_*=B(R^2)$ because $M(A_{\varepsilon})\cap \mathbf{E}_*\neq \{0\}.$
This causes $A_{\varepsilon}$ to be outside of $J(A, E_*).$ Therefore $A$ is not an inner point of $J(A, E_*)$ because $\lim\limits_{\varepsilon\rightarrow 0}A_\varepsilon=A,$ that is, $J(A,\mathbf{E}_*)$ is nontrivial.(The general result is given in Theorem $4.3$.)

{\bf Lemma $4.1$.}\quad{\it Suppose that $X$ is a double-split and non-zero operator in $B(E, F)$, say that $X^+$ is any of generalized inverses of $X$, and let $$\mathbf {E}_X=\left\{P_{N(X^+)}^{R(X)} T P_{N(X)}^{R(X^+)}: \forall T\in B(E, F)\right\}.$$ Then
  $$\mathbf{M}(X)=\left\{P^{N(X^+)}_{R(X)}T+P^{R(X)}_{N(X^+)}TP^{N(X)}_{R(X^+)}:\forall T\in B(E ,F)\right \},\eqno(4.1)$$
  $$\mathbf {E}_X \oplus \mathbf{M}(X)=B(E, F),$$
 and
  $$\mathbf {E}_X= \left \{T\in B(E,F): R(T)\subset N(X^+) \, \, \mbox{and} \,\, N(T)\supset R(X^+)\right\}.\eqno(4.2) $$}

\textbf{Proof :}\quad We first present the identity as follows,
$$ T=P^{N(X^+)}_{R(X)}T+P^{R(X)}_{N(X^+)}TP^{N(X)}_{R(X^+)}+ P_{N(X^+)}^{R(X)} T P_{N(X)}^{R(X^+)} \, \, \mbox{for all} \, \, T\in B(E,F).\eqno(4.3)$$
This is straightforward, but it is a key to proving the lemma.
By the definition of $\mathbf{M}(X)$,
$$P_{N(X^+)}^{R(X)}TP_{N(X)}^{R(X^+)}=0 \, \, \forall \, \, T\in M(X),$$ and by $(4.3)$, $$ T= P^{N(X^+)}_{R(X)}T+P^{R(X)}_{N(X^+)}TP^{N(X)}_{R(X^+)} \, \, \forall \, \, T\in \mathbf{M}(X).$$
 Thus, $$\{ P^{N(X^+)}_{R(X)}T+P^{R(X)}_{N(X^+)}TP^{N(X)}_{R(X^+)} : \, \forall T\in B(E, F)\}\supset M(X).$$ Conversely, for each $W\in B(E, F),$   $$(P^{N(X^+)}_{R(X)}W+P^{R(X)}_{N(X^+)}WP^{N(X)}_{R(X^+)})N(X)=P^{N(X^+)}_{R(X)}WN(X)\subset R(X) ,$$ therefore,
 $$\{ P^{N(X^+)}_{R(X)}T+P^{R(X)}_{N(X^+)}TP^{N(X)}_{R(X^+)} : \, \forall T\in B(E, F)\} \subset M(X).$$ Combining these two results , we confirm the equality $(4.1)$. It follows that $$\mathbf {E}_X \oplus \mathbf{M}(X)=B(E, F)$$ from the equalities $(4.3)$ and $(4.1).$ Now, remaining proof of the Lemma $4.1$ is verification of $(4.2).$ Obviously,
$$ R(P_{N(X^+)}^{R(X)} T P_{N(X)}^{R(X^+)})\subset N(X^+) \,\, \mbox{and} \,\, N(P_{N(X^+)}^{R(X)} T P_{N(X)}^{R(X^+)})\supset R(X^+),$$
these show that $\left \{T\in B(E,F): R(T)\subset N(X^+) \, \, \mbox{and} \,\, N(T)\supset R(X^+)\right\}\subset \mathbf{E}_X.$ Conversely, by $(4.3)$
$$ T=P_{N(X^+)}^{R(X)} T P_{N(X)}^{R(X^+)} \, \, \mbox{for} \, \, T\in  \left \{T\in B(E,F): R(T)\subset N(X^+) \, \, \mbox{and} \,\, N(T)\supset R(X^+)\right\},$$ this shows that  $\left \{T\in B(E,F): R(T)\subset N(X^+) \, \, \mbox{and} \,\, N(T)\subset R(X^+)\right\} \supset \mathbf{E}_X.$
Combining the above two results, we affirm the equality $(4.2)$. The proof ends. \quad$\Box$

{\bf Theorem 4.1.}\quad{\it Assume that $A$ is a non-zero, double-split operator in $B(E,F),$ then $\mathbf{\cal{F}}$ is smooth and integrable at $A.$}

{\bf Proof :}\quad  Let  $ \mathbf{V}_1^A=\{T\in B(E,F): \|(T-A)A^+\|<1\}$ and $\mathbf{V}_{1,A}=\mathbf{V}_1^A\cap \mathbf{V}(A^+, A)$.  Inspired by the conditions $(ii)$ in Theorem $1.1$, we define the following set  $\mathbf{S}:$
$$\mathbf{S}=\{T\in \mathbf{V}_{1,A}(A^+,A): A^+C_A^{-1}( A^+, T)\in GI(T)\},$$ where $GI(T)$ denotes the set of all generalized inverses of $T$, and $\mathbf{V}_{1,A}(A^+,A)$ is an open set containing $A.$ Evidently, for any $T\in \mathbf{S}$ , $A^+C_A^{-1}( A^+, T)= D_A^{-1}( A^+, T)A^+$ (see $(1)$ in the proof of theorem $1.1$ in the appendix) is the generalized inverse of $T$, denoted as $T^+, $ satisfying the conditions $N(T^+)=N(A^+)$ and $R(T^+)=R(A^+).$  To prove this theorem, consider the following mapping $\mathbf{\Phi}_A:$
 $$\mathbf{\Phi }_A(X)=(X-A)P_{R(A^+)}^{N(A)} + C_A^{-1}(A^+, X)X \, \, \mbox{for} \, \, X\in \mathbf{V}_1^A. \quad \eqno (4.4)$$
\, We claim that $\mathbf{\Phi}_A$ from $ \mathbf{V}_1^A$ onto itself is a $c^{\infty}$ mapping and satisfies $\mathbf{\Phi}_A(A)=A.$ To this end we introduce the two equalities as follows, $C_A(A^+, X)AA^+=XA^+$ and $$(\mathbf{\Phi}_A(X)-A)A^+= (X-A)A^+. \quad \eqno (4.5)$$ Indeed, for $X$ in $\mathbf{V}_1^A$, $$C_A(A^+, X)AA^+=AA^+ + (X-A)A^+AA^+=AA^+ + (X-A)A^+=XA^+$$ and $$(\mathbf{\Phi}_A(X)-A)A^+=(X-A)A^+ + C_A^{-1}(A^+, X)XA^+-AA^+ =(X-A)A^+.$$
 Since $C_A^{-1}(A^+, X)$ is of $c^{\infty},$  $\mathbf{\Phi}_A$ from $ \mathbf{V}_1^A$ onto itself is also a $c^{\infty}$ mapping with $\mathbf{\Phi}_A(A)=A.$ Moreover, $\mathbf{\Phi}_A(T)$ has the inverse as follows,
$$\mathbf{\Phi}_A^{-1}(T)=TP_{R(A^+)}^{N(A)} + C_A(A^+, T)TP_{N(A)}^{R(A^+)} \, \, \mbox{for all} \, \, T\in \mathbf{V}_1^A. \quad \eqno (4.6) $$
In fact, by $(4.5)$\begin{eqnarray*}
&&(\mathbf{\Phi}_A^{-1} \circ \mathbf{\Phi}_A)(X)\\&& = \mathbf{\Phi}_A(X)P_{R(A^+)}^{N(A)}+C_A(A^+, \mathbf{\Phi}_A(X))\mathbf{\Phi}_A(X)P_{N(A)}^{R(A^+)}
\\ &&=\mathbf{\Phi}_A(X)P_{R(A^+)}^{N(A)}+C_A(A^+, X)\mathbf{\Phi}_A(X)P_{N(A)}^{R(A^+)}\\ &&= (X-A)P_{R(A^+)}^{N(A)}+C_A^{-1}(A^+, X)XP^{N(A)}_{R(A^+)}+       C_A(A^+, X)(X-A)P_{R(A^+)}^{N(A)}P_{N(A)}^{R(A^+)}\\ &&+C_A(A^+, X)C_A^{-1}(A^+, X)XP_{N(A)}^{R(A^+)}+ C_A^{-1}(A^+, X)XP^{N(A)}_{R(A^+)} \\&&=(X-A)P_{R(A^+)}^{N(A)}+XP^{R(A^+)}_{N(A)}+C_A^{-1}(A^+, X)XP^{N(A)}_{R(A^+)}\\&&=X-A+C_A^{-1}(A^+, X)XP^{N(A)}_{R(A^+)};\end{eqnarray*}
 the equality $C_A(A^+, X)A(=XA^+A)=XP^{N(A)}_{R(A^+)}$ leads to $$(\mathbf{\Phi}_A^{-1} \circ \mathbf{\Phi}_A)(X)=X \, \, \mbox{for any} \,\, X\in \mathbf{V}_1^A.$$ Similarly,  $$(\mathbf{\Phi}_A \circ \mathbf{\Phi}_A^{-1})(X)=X \, \, \mbox{for any} \,\, X\in \mathbf{V}_1^A.$$
Combining the above results, we conclude $\mathbf{\Phi}_A$ from $ \mathbf{V}_1^A$ onto itself is a smooth diffeomorphism with $\mathbf{\Phi}_A (A)=A.$  Next,we discuss the relationship between $\mathbf{S}$ and $\mathbf{M}(A).$ The equivalence of conditions $(ii)$ and $(vi)$ in Theorem $1.1$ produces the following inferences and results:
\begin{eqnarray*}
\mathbf{S}&=&\{T\in \mathbf{V}_{1,A}: R(T)\cap N(A^+)=\{0\} \}\\&=&\{T\in \mathbf{V}_{1,A}: C_A^{-1}(A^+, T)T N (A)\subset R(A)\}\\&=&\{ T\in \mathbf{V}_{1,A}: \mathbf{\Phi }_A(T)N(A)\subset R(A)\}\\&=&\{T\in \mathbf{V}_{1,A}: \mathbf{\Phi }_A(T)\in \mathbf{M}(A)\}, \end{eqnarray*}  so  $$\mathbf{S}= \{T\in \mathbf{V}_{1,A}: \mathbf{\Phi }_A(T)\in \mathbf{M}(A)\}. \quad \eqno (4.7)$$
 For simplicity, write $\mathbf{M}(A)\cap \mathbf{\Phi}_A(\mathbf{V}_{1,A})$ as $\mathbf{V}_0.$ Based on $(4.7)$, we can conclude $$\mathbf{\Phi}_A(S)= \mathbf{V}_0 \, \, \mbox{ wich is an open set in} \,\, \mathbf{M}(A). \quad \eqno (4.8)$$
In fact, $(4.7)$ means $\mathbf{\Phi}_A(S)\subset \mathbf{V}_0 .$ For each $ Z\in \mathbf{V}_0$ let $T=\mathbf{\Phi }_A^{-1}(Z),$ then clearly, $T\in \mathbf{V}_{1,A},$ and $\mathbf{\Phi}_A(T)=Z\in \mathbf{M}(A),$ so $T$ belongs to $\mathbf{S}.$  This shows $\mathbf{\Phi}_A(S)\supset\mathbf{V}_0(A).$ Therefore, the property $(4.8)$ is true.
In addition,  by $(4.1)$ and $(4.2)$, $\mathbf {E}_A \oplus \mathbf{M}(A)=B(E, F).$ In summary, we demonstrate that $\mathbf{S}$ is a submanifold of $B(E,F)$ with only one coordinate chart. To complete the theorem's proof, we now need only to demonstrate that $T_X\mathbf{S}=\mathbf{M}(X)$. As a demonstration, we will prove that $T_A\mathbf{S}=\mathbf{M}(A).$ Let's start by finding the derivative of $\mathbf{\Phi}_A^{-1} .$ In order to simply let
 $$ \Gamma _A (A^+,X ,\triangle X)=C_A^{-1}(A^+, X+\triangle X)-C_A^{-1}(A^+, X).$$ By $(4.6)$
 \begin{eqnarray*}&&\mathbf{\Phi}_A^{-1}(T+\triangle T)-\mathbf{\Phi}_A^{-1}(T)\\ &&=(T+ \triangle T)P_{R(A^+)}^{N(A)} + C_A(A^+, T+\triangle T)(T+\triangle T)P_{N(A)}^{R(A^+)}-TP_{R(A^+)}^{N(A)} - C_A(A^+, T)TP_{N(A)}^{R(A^+)}\\&&=
 \triangle TP_{R(A^+)}^{N(A)}+\Gamma _A (A^+,T ,\triangle T)TP_{N(A)}^{R(A^+)}+ C_A(A^+, T+\triangle T)\triangle TP_{N(A)}^{R(A^+)}\\&&=\triangle TP_{R(A^+)}^{N(A)}+ \Gamma _A (A^+,T ,\triangle T)TP_{N(A)}^{R(A^+)}+\Gamma _A (A^+,T +\triangle T)\triangle T P_{R(A^+)}^{N(A)}+C_A(A^+, T)\triangle TP_{N(A)}^{R(A^+)}\\ && =\triangle TP_{R(A^+)}^{N(A)}+ \Gamma _A (A^+,T ,\triangle T)TP_{N(A)}^{R(A^+)}+C_A(A^+, T)\triangle TP_{N(A)}^{R(A^+)}
 +0(\|\triangle T\|),\end{eqnarray*} where $0(\|\triangle T\|)=\Gamma _A (A^+,T +\triangle T)\triangle T P_{R(A^+)}^{N(A)}$ because $\|\Gamma _A (A^+,T ,\triangle T)\|\rightarrow 0$ as $\|\triangle T\|\rightarrow 0.$  Therefore, $$(\mathbf{\Phi}_A^{-1})'(T)\triangle T=\triangle TP_{R(A^+)}^{N(A)}+ \Gamma _A (A^+,T ,\triangle T)TP_{N(A)}^{R(A^+)}+C_A(A^+, T)\triangle TP_{N(A)}^{R(A^+)};$$ specially,
 $$(\mathbf{\Phi}_A^{-1})'(A)=\mathbf {I}$$
where $\mathbf {I}$ represents the identity on $B(E, F).$ Using a proof method similar to that in formula $(2.3),$ it can be shown that
  $$(\mathbf{\Phi}_A^{-1})'(A) T_A \mathbf{S}=\mathbf{M}(A) \,\, \mbox{for} \, \, X\in \mathbf{S},$$
  which leads to
$T_A \mathbf{S}=\mathbf{M}(A).$ Next, we are going to prove that $T_B \mathbf{S}= \mathbf{M}(B)$ for any $B \in \mathbf{S}.$
Let $\mathbf {V}_1^B=\{ T\in B(E, F): \|(T-B)B^+\|<1\},$ and
 $$\mathbf{\Phi }_B(X)=(X-B)P_{R(A^+)}^{N(B)} + C_B^{-1}(B^+, X)X \, \, \mbox{for} \, \, X\in \mathbf{V}_1^B.$$
Like $\mathbf{\Phi}_A^{-1}(T),$ we also have
$$\mathbf{\Phi}_B^{-1}(T)=TP_{R(A^+)}^{N(B)} + C_B(B^+, T)TP_{N(B)}^{R(A^+)} \, \, \mbox{for all} \, \, T\in \mathbf{V}_1^B,$$
and
$\mathbf{\Phi}_B$ from $ \mathbf{V}_1^B$ onto itself is a $c^{\infty}$ diffeomorphism with that $\mathbf{\Phi}(B)=B$ and
$(\mathbf{\Phi}_B^{-1})'(B)=\mathbf {I}.$ Let
$\mathbf{W}_{1,B}=\mathbf{V}_1^B\cap \mathbf{V}(B^+, B)\cap \mathbf{V}_{1,A}.$ Clearly,  $B$ belongs to $\mathbf{W}_{1,B}\cap \mathbf {S},$  and $ N(B^+)=N(A^+).$ Then
\begin{eqnarray*}
\mathbf{S\cap \mathbf{W}_{1,B}}&=&\{T\in \mathbf{W}_{1,B}: R(T)\cap N(A^+)= \{0\} \}\\&=&\{T\in \mathbf{W}_{1,B}: R(T)\cap N(B^+)= \{0\} \}\\&=&\{T\in \mathbf{W}_{1,B}: C_B^{-1}(B^+, T)T N (B)\subset R(B)\}\\&=&\{ T\in \mathbf{W}_{1,B}: \mathbf{\Phi }_B(T)N(B)\subset R(B)\}\\&=&\{T\in \mathbf{W}_{1,B}: \mathbf{\Phi }_B(T)\in \mathbf{M}(B)\},\end{eqnarray*} so
$$ \mathbf{\Phi }_B^{-1}(\mathbf{S \cap \mathbf{W}_{1,B} })=\mathbf{M}(B)\cap  \mathbf{\Phi }_B( \mathbf{W}_{1,B}). $$
 The method used to prove $T_A \mathbf{S}=\mathbf{M}(A)$ can also derive
 $$T_B \mathbf{S}=(\mathbf{\Phi}_B^{-1})'(B) T_B \mathbf{S}=\mathbf{M}(B) \, \mbox{for any} \, B\in \mathbf{S}.$$
The theorem is proved. $\Box$

{\bf Note:}  The operator $A$ in the theorem can be highly singular in which case $\mbox{dim} N(A)=\mbox{dim} R(A)=\mbox{codim} R(A)=\infty.$  The theorems $4.1$ and $4.3$ show that $\mathbf{S}$ connects it into a class, especially its co-tailed set is non-trivial. These are significant new results in operator topology, as well as in global analysis.

{\bf Theorem 4.2}\quad{\it Let $\mathbf{S}$ denote any of the following classes:
$$F_k,  \Phi_{m,n}, \Phi_{m,\infty} \, \, \mbox{and} \, \, \Phi _{\infty, n},$$  $k=1, 2,\cdots ,$ and $m,n=0, 1, \cdots . $ Then  $\mathbf{S}$ is a smooth submanifold in $B(E,F)$, and is tangent to $\mathbf{M}(Z)$ at any $Z$ in $\mathbf{S}$.}

 {\bf Proof}\quad   It is well known that any $Z\in \mathbf{S}$ is double-split, say that $Z^+$ is a generalized inverse of $Z$. Write $ \mathbf{V}_1^Z=\{X\in B(E, F):\|(T-Z)Z^+\|<1\}$ and $\mathbf{V}_{1,Z}=\mathbf{V}_1^Z \cap \mathbf{V}(Z^+, Z).$

 By Theorems $1.1$ and $1.2,$
 $$\{T\in \mathbf{V}_{1,Z}: Z^+C_Z^{-1}( Z^+, T)\in GI(T)\}=\mathbf{S}\cap \mathbf{V}_{1,Z} \, \, \mbox{for any Z} \, \, \in \mathbf{S}. $$
 Similar to $(4.4)$ and $(4.6)$, we also have that
  $$\mathbf{\Phi}_Z(T)=(T-Z)P_{R(Z^+)}^{N(Z)}+C^{-1}_Z(Z^+,T)T \, \, \mbox{for any} \, \, T\in  \mathbf{V}_1^Z .\quad \eqno (4.9),$$ and
 $$\mathbf{\Phi}_Z^{-1}(T)=TP_{R(Z^+)}^{N(Z)} + C_Z(Z^+, T)TP_{N(Z)}^{R(Z^+)} \, \, \mbox{for any T} \, \, \mathbf{V}_1^Z . \quad \eqno (4.10)$$
According to Theorem $4.1,$ $\mathbf{S}\cap \mathbf{V}_{1,Z} $ is a smooth submanifold tangent to  $\mathbf{M}(X)$ at any $X$ of it, and its atlas has only one coordinate chart $(\mathbf{\Phi}_Z|_{\mathbf{S}}, \mathbf{S}\cap \mathbf{V}_{1,Z},\mathbf{M}(Z)).$
 To complete the proof of the theorem, we now need only prove that the atlas $\{(\mathbf{\Phi}_Z, \mathbf{S}\cap \mathbf{V}_{1,Z}, \mathbf{M}(Z))\}_{Z\in \mathbf{S}}$ is compatible.Suppose that for $A$ and $Z$ in $\mathbf{S},$  $\mathbf{S} \cap \mathbf{V}_{1,A} \cap \mathbf{V}_{1,Z}\neq\emptyset $.
Let $\mathbf{O}_Z=\mathbf{\Phi}_Z(\mathbf{S} \cap \mathbf{V}_{1,Z} \cap \mathbf{V}_{1,A})$, and $\mathbf{O}_A=\mathbf{\Phi}_A(\mathbf{S} \cap \mathbf{V}_{1,Z} \cap \mathbf{V}_{1,A})$. Obviously, $\mathbf{O}_Z$ and $\mathbf{O}_A$ are open sets in $\mathbf{M}(Z)$ and $\mathbf{M}(A)$, respectively. Clearly, $$ \mathbf{\Phi}_Z\circ \mathbf{\Phi}_A^{-1} \, \, \mbox{from $\mathbf{O}_A$ onto $\mathbf{O}_Z $ is  $c^{\infty}$ diffeomorphism }.$$
 This is immediate from the equalities $(4.9)$ and $(4.10)$.This shows that the atlas $\{(\mathbf{\Phi}_X|_{\mathbf{S}}, \mathbf{S}\cap \mathbf{V}_{1,X}, \mathbf{M}(X))\}_{X\in \mathbf{S}}$ is compatible. The proof of the theorem ends. $\Box$

 {\bf Theorem 4.3.}\quad{\it Assume that $A_0$ belongs to $B^+(E, F)$, and that $\dim N(A_0)$ and $\dim N(A_0^+)$ are positive numbers. Write $\mathbf {E}_{A_0}=\mathbf{E_*}.$ Then the co-tailed set of $\mathbf{\cal{F}}$ at $A_0$,  $J(A_0,\mathbf{E}_*)$  is non-trivial.}

 \textbf{Proof :}\quad   Let both  $e$ in $N(A_0)$ and $e^+$ in $ N(A_0^+)$ be non-zero elements, $[e]$ denote the one-dimensional space generated by $e$, and $N_0^*$  satisfy $N(A_0)=N_0^*\oplus [e].$  Define a non-zero operator $L$ in $B(E,F),$ satisfying that
$$N(L)=N_0^*, \, \, L|_{[e]}\in B^+([e],[e^+]),$$ and for each  $ x\in R(A_0^+),$ $$Lx = A_0x .$$
For arbitrary  positive number $\varepsilon,$ consider the operator as follows,
$$A_\varepsilon = A_0+\varepsilon L .$$ We claim that $$N(A_\varepsilon)=N_0^*.$$  For any $x \in E$ let
$x = x_* + x_1 + x_+$ where $x_* \in N_o^*, x_1\in [e]$ and $x_+ \in R(A_0^+).$ By computing directly,
\begin{eqnarray*} N(A_\varepsilon)&=&\{  x_* + x_1 + x_+ : A_\varepsilon (x_* + x_1 + x_+)=0\}\\&=& \{ x_* + x_1 + x_+ : A_0 x_+ + \varepsilon L x_1+ \varepsilon A_0x_+=0\}\\&=&\{ x_* + x_1 + x_+ : (1+\varepsilon )A_0x_+ +\varepsilon L x_1 = 0 \}\\&=&\{ x_* + x_1 + x_+ :  (1+\varepsilon )A_0x_+=0 \, \, {and} \, \, \varepsilon L x_1=0 \}\\ &=& \{ x_* : \forall x_* \in N_0^* \}.
\end {eqnarray*} This shows $N(A_\varepsilon)=N_0^*.$

Consequently, for any $ A \in M(A_0)$
$$AN(A_\varepsilon)=A N_0^*\subset R(A_0)\subset R(A_\varepsilon),$$
and so, $$M(A\varepsilon)\supset M(A_0), \, \, \rm {but} \, \, M(A_\varepsilon)\neq M(A_0) \quad \eqno (4.11)$$
because of $A_\varepsilon e$ not being in $R(A_0)$. Obviously,
  $$\parallel A_\varepsilon- A_0\parallel =\mid\varepsilon\mid  \parallel L \parallel ,$$ and hence,
$$\lim\limits_{\varepsilon\rightarrow 0}A_{\varepsilon}=A_0. \quad \eqno (4.12)$$
To conclude the proof, we next want to show that $A_\varepsilon$ does not belong to $J(A_0,E_*), $  Otherwise, suppose that $A_{\varepsilon}$ is in $J(A_0,E_*). $
Then
 $$ M(A_{\varepsilon})\oplus E_*= M(A_0)\oplus E_*= \rm{B(E,F)}.\quad \eqno (4.13)$$
Consequently, for any $X$ in $ B(E,F)$,
$$ X=P_{M(A_{\varepsilon})}^{E_*}X+P^{M(A_{\varepsilon})}_{E_*}X = P_{M(A_0)}^{E_*}X+P^{M(A_0)}_{E_*}X $$
and
$$( P_{M(A_{\varepsilon})}^{E_*}X-P_{M(A_0)}^{E_*}X )+(P^{M(A_{\varepsilon})}_{E_*}X- P^{M(A_0)}_{E_*}X)=0.$$
Since that $M(A_{\varepsilon})\supset M(A_0)$ and $M(A_{\varepsilon})\oplus E_*= \rm{B(E,F)}$( see $(4.11)$ and $(4.13)$ ),
 $$P_{M(A_{\varepsilon})}^{E_*}X=P_{M(A_0)}^{E_*}X \, \, \mbox{for  any  X  in  B(E,F)}.$$
This means $M(A_{\varepsilon})=M(A_0),$ which is contradictory to $(4.11)$. Therefor $A_{\varepsilon}$ does not belong to $J(A_0,E_*). $
By $(4.11)$ and $(4.12)$ we conclude that $A_0$ is not the inner point of $J(A_0, E_*)$, meaning that $J(A_0, E_*)$ is non-trivial.\quad$\Box$.

Let $U$ be an open set containing point $x_0,$ $\mathbf{S}$ a $c^1$ submanifold in $U,$ and $f$ from $U$ to $(\infty, -\infty)$ a $c^1$ map, we have

 {\bf Theorem 4.4}\quad{\it If $x_0\in U\cap \mathbf{S}$ is an extreme point of $f$ under the constrain of submanifold $\mathbf{S},$ then
$ N(f'(x_0))\supset T_{x_0}\mathbf{S},$ that is, $ f'(x_0)T_{x_0}\mathbf{S}=0.$}

\textbf {Proof :}\quad According to definition of $c^1$ submanifold $\mathbf{S}$ in $U,$ there exists a split subspace $E_0$ in $E$, a neighborhood $U_0$ at $x_0$ and a $c^1$ diffeomorphism $\varphi : U_0 \rightarrow \varphi (U_0),$ such that $\varphi(\mathbf{S}\cap U_0)$ is an open set in $E_0.$ Without loss of generality, still write $U\cap U_0$ as $U_0.$ We claim that
$$ E_0=\varphi'(x_0)^{-1}T_{x_0}\mathbf{S}. \quad \eqno (4.14)$$
By the definition of the tangent space of $\mathbf{S}$ at $x_0,$  $ T_{x_0} \mathbf{S}=\{\dot{c_{(\delta, x_0}}(0): \, \, \forall c^1 \, \, \mbox{curve} \, \, c_{\delta,x } (t)\},$ where $c_{\delta,x_0}(t)$ is as the same as in the proof of Theorem $2.1.$ Since $\varphi (S\cap U_0)$ is an open set in $E_0,$  $$\frac{d}{dt}(\varphi \circ c_{\delta, x})\mid_{t=0}=\varphi'(x)\dot{c_{\delta, x}}(0)\in E_0 \, \, \mbox{for any $c^1$ curve} \, \, c_{\delta, x}(t),$$ therefore $\varphi'(x)T_x \mathbf{S}\subset E_0.$ Conversely, let $\delta$ be a sufficiently small positive number $\delta$, and satisfy that for any $t\in (-\delta, \delta )$, $\varphi (x)+t e \in \varphi (S\cap U_0)$. Directly,
$$\dot{c_{\delta, x}}(0)=(\varphi ^{-1})'(\varphi(x))e=\varphi'(x)^{-1}e \, \, \mbox{for \, all} \, \, e\in E_0,$$
therefore $\varphi'(x)T_x \mathbf{S}\supset E_0.$ Combining the two results above we prove the equality $(4.14).$ Now go back and prove the theorem.We convert the extreme point $x_0$ of $f\mid_{\mathbf{S}\cap U_0}$ into extreme point $\varphi (x_0)$ of unconstrained functional $f\circ \varphi ^{-1}$ defined on the open set $\varphi (\mathbf{S}\cap U_0)$
in $E_0.$ Therefore $(f \circ \varphi ^{-1})'(\varphi(x_0))=0$ in $B(E_0, R).$ Evidently,
$$f'(x_0)\cdot(\varphi ^{-1})'(\varphi(x_0)e=f'(x_0)\varphi'(x_0)^{-1}e = 0 \, \, \mbox{for} \, \, e \in E_0.$$ Let $\sigma = \varphi'(x_0)^{-1}e ,$ then by $(4.14)$ $$f'(x_0)\sigma =0 \, \, \mbox{ for  any } \, \, \sigma \in T_{x_0}\mathbf{S}.$$ The proof ends.$\Box$

{\bf Note:} The generalized transversality theorem (see [Ma5]) is a global implicit function theorem, particularly featuring a tangent space formula. Consequently, Theorem $4.4$ finds extensive applications, as demonstrated in [Ma7].

\begin{center}\textbf{Appendix}
\end{center}

 \textbf{$1$ \, The proof of Theorem $1.1.$}

 Before proving Theorem $1.1,$ we give the following equalities for $T \in V (A, A^+):$
$$B=A^+C^{-1}_A(A^+,T)=D^{-1}_A(A^+,T)A^+,  BTB=B, \, \mbox{and} \, C^{-1}_A(A^+,T)TP_{R(A^+)}^{N(A)}=A. \eqno (1)$$ Indeed,
\begin{eqnarray*}
&& A^+C^{-1}_A(A^+,T)- D^{-1}_A(A^+,T)A^+ \\ && = D^{-1}_A(A^+,T)(D_A(A^+,T)A^+ - A^+ C_A(A^+,T))C^{-1}_A(A^+,T)\\
 &&=D^{-1}_A(A^+,T)(A^+T - A^+T)C^{-1}_A(A^+,T)=0;\end{eqnarray*}
 \begin{eqnarray*}
B T B &=&A^+ C^{-1}_A(A^+,T)TA^+C^{-1}_A(A^+,T)\\ &=&A^+C^{-1}_A(A^+,T)(C_A(A^+,T)- P_{N(A^+)}^{R(A)})C^{-1}_A(A^+,T) \\&=&A^+C^{-1}_A(A^+,T)\, \, \mbox{because $C_A(A^+,T)P_{N(A^+)}^{R(A)}=P_{N(A^+)}^{R(A)}$}\\ &=&B, \, \, \mbox{(also refer to [N-C])};\end{eqnarray*}  the third equality in $(1)$ follows from the obvious equality $ C_A(A^+,T)A=T P_{R(A^+)}^{N(A)}.$ Now we begin to prove the theorem.
$$\mathbf{(vii)}\Leftrightarrow \mathbf{(ii)}.$$
 Evidently,
 \begin{eqnarray*}
T-T B T &=& (C_A(A^+, T)-TA^+)C^{-1}_A(A^+, T)T\\&=&
P_{N(A^+)}^{R(A)}C^{-1}(A^+,T)T \, \quad \mbox{for} \, \, T\in V(A,A^+). \end{eqnarray*}
From this, we can deduce that $(vii)\Leftrightarrow(ii).$
$$\mathbf{(vi)} \Leftrightarrow \mathbf{(vii)}.$$
By the third equality in $(1),$
$$\begin{array}{rllr}
C^{-1}_A(A^+,T)Th&=&C^{-1}_A(A^+,T)TA^+Ah+C^{-1}_A(A^+,T)T(I_E-A^+A)h\\
&=&Ah+C^{-1}_A(A^+,T)T(I_E-A^+A)h \quad\forall h\in E.\end{array}$$
This leads to the equivalence of $(vi)$ and $(vii).$
$$\mathbf{(v)}\Leftrightarrow \mathbf{(vi)}$$
To verify that $(vi) \Rightarrow (v),$ assume that $(vi)$ is true, in other words, for each $h \in N (A)$ there exists a $g\in R(A^+)$ such that $Th=C_A(A^+,T)Ag.$ Obviously, $C_A(A^+,T)Ag=Tg ,$ so $h-g$ belongs to $N(T),$ which leads to $(I_E-A^+A)(h-g)=h \, \, \mbox{for any h} \in N(A).$
Now we can conclude $(vi) \Rightarrow (v).$
To verify that $(v)\Rightarrow (iv) $ we assume that $(v)$ is true, and for each $h\in N(A),$ let $g\in N(T)$ satisfy $h=(I_E-A^+A)g$. Then
 $C^{-1}_A(A^+,T)Th=-C^{-1}_A(A^+,T)TA^+Ag=-Ag\in R(A).$ This indicates that the condition $(vi)$
is valid, therefore, $(v)$ $\Leftrightarrow$ $(vi).$
$$\mathbf{(i)}\Leftrightarrow \mathbf{(ii)}.$$
To show that  $(ii)\Rightarrow (i)$ assume that the condition $(ii)$ is true. Then $B$ is a generalized inverse of $T,$ and satisfies $N(B)=N(A^+)$. Therefore, the condition $(i)$ holds.
To show that $(i)\Rightarrow (ii)$ assume that for $T\in V(A, A^+)$, $R(T)\cap  N(A^+)= \{0\}.$ It has been pointed out in the proof of
$(vii)\Leftrightarrow(ii)$ that  $ T - T B T =  P_{N(A^+)}^{R(A)}C^{-1}(A^+,T)T,$  thus $R( T - T B T)\cap N(A^+)\subset R(T)\cap N(A^+)=\{0\}.$ This shows $T-T B T=0,$ and hence, $B$ is the generalized inverse of $T,$ satisfying that $R(B)=R(A^+)$ and $N(B)=N(A^+).$
 $$\mathbf{(i)}\Leftrightarrow \mathbf{(iii)}.$$
Obviously,
$(iii)$ $\Rightarrow$ $(i).$  We only need to verify the converse statement. Assuming that the condition $(i)$ is true,from $(i)$ $\Leftrightarrow$ $(ii)$, we can deduce that $B$ is the generalized inverse of $ T.$  Therefore,
  $N (B) =N (A ^ +)$ and $ F=R (T) \oplus N (A ^ +).$  This indicates that $(i)$ $\Rightarrow$ $(iii).$ The proof of $(i)$ $\Leftrightarrow$ $(iii)$ ends.
$$\mathbf{(i)}\Leftrightarrow \mathbf{(iv)}.$$
Assume that the condition $(i)$ is true, then $E=N(T)\oplus R(B)=N(T)\oplus R(A^+)$ because $(i)\Leftrightarrow(ii).$ This shows that $(i)\Rightarrow (iv)$. Conversely,
assuming that the condition $(iv)$ is true, $N(A)=(I_E-A^+A)E=(I_E-A^+A)(N(T)\oplus R(A^+))=(I_E-A^+A)N(T).$ This shows that the conditions $(v)$ is valid, so that $(iv)$ $\Rightarrow $ $(v)$.

The following equivalence relationships have been previously demonstrated:
$$ (v) \Leftrightarrow (vi), (vi) \Leftrightarrow (vii), (vii) \Leftrightarrow (ii), \, \mbox{and} \, (ii) \Leftrightarrow (i).$$
 From these we can deduce that $(iv)\Rightarrow (i).$ So $(i) \, \, \Leftrightarrow \, \, (iv).$

Summarizing the following equivalence relationships shown above: $(i) \, \, \Leftrightarrow \, \, (ii)$, $(i) \, \, \Leftrightarrow \, \, (iii)$, $(i) \, \, \Leftrightarrow \, \, (iv),$
 $(ii) \, \, \Leftrightarrow \, \, (vi),$
 $(v) \, \, \Leftrightarrow \, \, (vi),$ and $(vi) \, \, \Leftrightarrow \, \, (vii),$ we assert that
 the theorem $1.1$ is true.

\textbf{$ 2$ \, The proof of Theorem $1.2$}

Let $A\in B^+(E, F),$  $A^+$ be any one of generalized inverses of $A,$ and  $B=A^+C^{-1}_A(A^+,T) $ for $T\in V(A, A^+)$.
For each $T\in  V(A, A^+)$, operators $B$ and $T$  produce the following two interesting projections
 $P_1$ and $P_2$:
$$P_1=BT \, \, \, \mbox{and} \, \, \, P_2=TB.$$
 By the second equality in $(1),$
$$P_1^2=BTBT=BT=P_1  \, \, \mbox{and} \, \, P_2^2 =TBTB=TB=P_2.$$
 This indicates that $P_1$ and $P_2$ are projections on $E$ and $F$, respectively. Furthermore,
 $$ P_1= P^{N(A^+T)}_{R(A^+)} \, \, \mbox{and} \, \,  P_2= P_{R(TA^+)}^{N(A^+)} \, \, \, \forall T\in V(A, A^+) .\quad \eqno (2) $$
Indeed, by the first equality in $(1)$,
$$N(P_1)=N(BT)=N(D_A^{-1}(A^+,T)A^+T)=N(A^+T) \, \, \mbox{for} \, \, T\in V(A, A^+);$$
by the third equality in $(1),$ $$P_1A ^ + =A ^ + C ^ {-1} _A (A ^ +, T) TA ^ + =A ^ + AA ^ + =A ^ + =A ^ +,$$ indicating that $R (A ^ +) \subset  R (P_1),$ and obviously, $R(P_1)\subset R(A^+),$ therefore,
  $R (P_1) = R (A ^ +).$  Combined with the above conclusions for $N(P_1)$ and $R(R_1)$ we prove the first formula in $(2).$ Similarly,
$$R(P_2)= R(TA^+C_A^{-1}(A^+,T))=R(TA^+);$$
 due to the second equality in $(1),$ for each $e\in N(P_2),$ $$B e=B T B e =B P_2 e =0,$$ and obviously, $N(P_2)\subseteq N(A^+),$
therefore, $N(P_2) =N(A^+).$ Now the proof of the two formulas in $(2)$ is complete.
  These two equalities are equivalent to
 $$  E=R(A^+)\oplus N(A^+T) \, \, \mbox{and} \, \,  F=R(TA^+)\oplus N(A^+). \quad \eqno(3)$$
Next, we began to prove Theorem $1.2.$

  $(a)$  We begin to prove Theorem $1.2$ in the case where  $\dim R(A)<\infty .$ The following property about $T \in V(A^+,A)$ will be needed: $$R(T)\cap N(A^+)=\{o\} \, \, \mbox{if and only if} \, \, \dim R(T)=\dim R(A)<\infty .$$ Now let us prove this. In fact, if $R (T) \cap N (A ^ +) = \{0\}$ for $T\in V (A ^ +, A),$ then by the equivalence of conditions $(i)$ and $(iii)$ in Theorem $1.1,$ $F= R (T)\oplus N (A ^ +)= R(A)\oplus N(A^+).$ Thus, $\mbox{dim} R (T) =\mbox{dim} R (A)< \infty .$ Proving the reverse is not so simple. For $T\in B^+(E, F)$ let $T^+$ be a generalized inverse of $T$. We introduce the following subspace $E_*$ in $R(T^+)$: $$E_*=\{e\in R(T^+): Te\in N(A^+)\}, \, \mbox{ensuring that} \, N(A^+T)= N(T)\oplus E_*.$$ By the latter equality in $(3)$, $ F=R(TA^+)\oplus N(A^+)= R(A)\oplus N(A^+),$ and so, $\mbox{dim} R(TA^+)=\mbox{dim} R(T)< \infty.$ It is time to prove the reverse. If
  $ \mbox{dim} R(T)=\mbox{dim} R(A)<\infty ,$ then the first equation in $(3)$ imply  $R(T)=R(TA^+)\oplus TE_*,$ which consequently leads to $\mbox{dim} TE_*=0.$ Therefore, since $E_*\subset R(T^+),$ $\mbox{dim} E_*=0.$ This shows that for $T\in V(A, A^+),$ $R(T)\cap N(A^+)=\{0\}.$ Combining the above two results, we prove that the theorem $1.2$ holds when the operator $A$ belongs to any of the following classes : $F_k, k=1, 2, \cdots.$

$(b)$ Assume that $A\in B^+(E, F)$ and $\dim N(A)<\infty.$ The following property about $T\in V(A, A^+)$ will be needed : $R(T)\cap N(A^+)=\{0\} $ is equivalent to the following three conditions:  $$ T\in V^+(A, A^+), \, \mbox{dim} N(T)=\mbox{dim} N(A)< \infty , \, \, \mbox{and}\, \mbox{codim} R(T) = \mbox{codim} R(A),$$
where $ V^+(A, A^+)=B^+(E, F)\cap V(A, A^+).$ Now let us prove this. In fact, If  $R(T)\cap N(A^+)=\{0\} \, \mbox{for any} \, T\in V(A^+,A),$ then by the equivalence of the conditions  $(i)$ and $(iv)$ in Theorem $1.1$, $ E=N(T)\oplus R(A^+)=N(A)\oplus R(A^+),$ and hence, $\mbox{dim} N(T)= \mbox{dim} N(A)<\infty $ ;
by the equivalence of the conditions  $(i)$ and $(iii)$ in Theorem $1.1,$  $F=R(T)\oplus N(A^+)=R(A)\oplus N(A^+),$ and hence, $\mbox{codim} R(T)=\mbox{codim} R(A);$ the equivalence of the conditions $(i)$ and $(ii)$ in the theorem $1.1$ leads to $ T \in V^+(A, A^+).$ Summarizing the above three conclusions we assert that if $T$ is in $V(A,A^+)$, and $R (T) \cap N (A ^+) =\{ 0\},$ then $T$ belongs to $V^+(A,A^+)$ and satisfies that $\mbox{dim} N (T) = \mbox{dim} N(A)<\infty $ and $\mbox{codim} R(T) =\mbox{codim} R(A).$ Conversely, assume that $ T\in V^+(A, A^+), \, \,\mbox{dim} N(T)=\mbox{dim} N(A)<\infty $ and $\mbox{codim}R(T)=\mbox{codim}R(A),$ we want to show that $ T\in V(A^+, A)$ and $R(T)\cap N(A^+)=\{0\}.$ By the equivalence of the conditions $(i)$ and $(ii)$ in Theorem $1.1,$ the condition $(i)$ includes that $T$ has a generalized inverse $B$ satisfying $N(B)=N(A^+)$, hence $\mbox{codim}R(T)=\mbox{codim }R(A)$.  Therefore, we only need to show that if for $T\in V^+(A, A^+),$ $\mbox{dim} N(T) =\mbox{dim} N(A) <\infty ,$ then $ T\in V^+(A, A^+)$ and $R(T)\cap N(A^+)=\{0\}.$ For $T\in V^+(A, A^+),$ let $T^+$ be a generalized inverse of $T.$ By the two equalities,$N(A^+T) =E_* \oplus N(T)$ and the previous equality in $(3)$,
$$ E=R(A^+)\oplus N(A^+T)=R(A^+)\oplus N(A),$$ so that $$\mbox{dim}N(A^+T)=\mbox{dim}N(A)=\mbox{dim}N(T)<\infty.$$ therefore, $\mbox{dim} E_*=0,$ which results in that $T \in V (A , A^+)$ and $R (T) \cap N (A ^ +) = \{0\}.$ Now, we prove that the theorem $1.2$ holds when $A$ belongs to any of the following classes: $\Phi_{m,n}$ and $\Phi_{m, \infty}, m =0, 1, \cdots.$

$(c)$ Assuming that $A\in B^+(E,F)$ and $\mbox{codim}R(A)<\infty,$ we deed to verify the following property : $T\in V(A, A^+)$ and
$R(T)\cap N(A^+)=\{0\}$ are equivalent to the following three conditions: $$ T\in V^+(A, A^+), \, \,  \mbox{codim} R(A) = \mbox{codim} R(T) \,< \infty , \, \, \mbox{and} \, \,  \mbox{dim} N(A)= \mbox{dim} N (T).$$
Assuming that $R(T)\cap N(A^+)=\{0\}$ and $T\in V(A, A^+)$, based on the three mutually equivalent conditions $(i),$ $(iii),$ and $(iv)$ in Theorem $1.1$, we deduce that for any $T\in V^+(A, A^+)$, $ E=N(T)\oplus R(A^+)=N(A)\oplus R(A^+)$ and $F=R(T)\oplus N(A^+)=R(A)\oplus N(A^+).$
This demonstrates that if  $R(T)\cap N(A^+)=\{0\}$ and $T\in V(A, A^+),$ then $T$ belongs to $V^+(A, A^+),$ furthermore, $\mbox{codim} R(A)$ = $\mbox{codim} R(T)$ $<\infty $ and $\mbox{dim} N(A)=\mbox{dim} N(T).$
Conversely,  assume that $T\in V^+(A, A^+)$,  $\mbox{codim} R(A)$ = $\mbox{codim} R(T)< \infty $ and $\mbox{dim} N(A)=\mbox{dim} N(T).$  We want to prove that  $T\in V(A, A^+)$ and $R(T)\cap N(A^+)=\{0\}.$ Let $T^+$ be a generalized inverse of $T.$ According to the two equalities $ N(A^+T)=E_* \oplus N(T)$ and the previous equality in $(3)$,
 $$R(T)=R(TA^+)\oplus TE_*, \, \, \mbox{for} \, \, T\in V^+(A, A^+),$$ and hence,
$$ F=R(TA^+)\oplus TE_*\oplus N(T^+).$$  In addition, by the latter formula in $(3)$
$$ F= R(T)\oplus N(T^+)= R(TA^+)\oplus TE_T^*\oplus N(T^+)=R(TA^+)\oplus N(A^+),$$ and hence $$ \mbox{dim} N(A^+)=\mbox{dim} TE_*+\mbox{dim} N(T^+).$$
Consequently, the assumption, $\mbox{codim} R(A) = \mbox{codim} R(T)<\infty$ i.e., $\mbox{dim} N(A^+)= \mbox{dim} N(T^+) < \infty$ leads to $\mbox{dim} TE_*=0 .$  So, $T\in V(A, A^+)$ and $R(T)\cap N(A^+)=\{0\}.$ Now, we prove that the theorem $1.2$ holds when $A$ belongs to any of the following classes : $\Phi_{m,n}$ and $\Phi_{\infty, n} ,m.n=0, 1, \cdots.$ Summarizing the three conclusions $(a)$ $(b)$ and $(c)$ above, we prove Theorem $1.2$.

\textbf{$3$\quad The Proof of Theorem $1.3$ }
Let $T_x$ be an operator valued map from a topological space $X$ into $B(E,F)$, and be continuous at $x_0 \in X$.
 Assume that
 $$R(T_x)\cap N(T^+_0)=\{0\} \, \, \mbox{for all} \, \, x\in U_0,$$ where $U_0$ is a neighborhood at $x_0$, $T_0=T_{x_0}$, and $T^+_0$ is a generalized inverse of $T_0$. Let $T^{\oplus}$ be any other generalized inverse of $T_0$, $\delta = min \{\|T^+_0\|^{-1},\|T^+_0T_0T^\oplus\|^{-1}\|\},$
 and $V_\delta=\{T\in B(E,F):\|T-T_0\|<\delta\}.$
 According to the continuity of $T_x$ at $x_0$,
There exists in $U_0$ a neighborhood $U_1$ at $x_0$ such that for $x\in U_1,$ $T_x$ belongs to $V_\delta .$
 For simplicity, still write $U_1$ as $U_0$, and then, we have that for all $x\in U_0, R(T_x)\cap N(T^+_0)=\{0\}$ and $T_x\in V_\delta .$
We claim that $T^\oplus $ also meets that $$R(T_x)\cap N(T^\oplus )=\{0\}\quad\forall x\in U_0.$$ Write $B=T^+_0T_0T^\oplus.$
Obviously,$BT_0B=T^+_0(T_0T^\oplus T_0)T^+_0T_0T^\oplus=T^+_0(T_0T^+_0T_0)T^\oplus=T^+_0T_0T^\oplus=B;$
$T_0BT_0=T_0T^+_0(T_0T^\oplus T_0)=T_0T^+_0T_0=T_0$.
These show that $B$ is also a generalized inverse of $T_0$.Further,$R(B)=R(T_0^+)$ and $N(B)=N(T^\oplus).$ Indeed, $R(B)=R(BT_0)=R(T_0^+T_0T^\oplus T_0)=R(T_0^+T_0)=R(T_0^+)\mbox{and} N(B)=N(T_0B)=N(T_0T_0^+T_0T^\oplus)=N(T_0T^\oplus)=N(T^\oplus).$
Taking $A$ and $A^+$ in Theorem $1.1$ as $T_0$ and $T^+_0$ respectively,
The equivalence of conditions $(i)$ and $(iv)$ in the correspondence theorem leads to the result as follows,
$R(T^+_0)\oplus N(T_x)=E \, \forall x\in
U_0.$ That is, $$R(B)\oplus N(T_x)=E \quad\forall x\in U_0.$$
Similarly, instead of $A$ and $A^+$ in Theorem $ 1.1$ by $T_0$
and $B$, respectively. The equivalence of conditions $(iv)$ and $(i)$ in the corresponding theorem shows that $R(T_x) \cap N(B)=\{0\}.$ Therefore,
$$R(T_x)\cap N(T^\oplus )=\{0\} \, \, \mbox{for all} \, \, x\in
U_0.$$ The proof ends.

\textbf{$4$ \quad The Proof of Theorem $1.4$}

Let $T_0=T_{x_0},$ and $T^+_0$ be any generalized inverse of $T_0.$ Suppose that $x_0$ is a local fine point of $T_x$, and $U_0$ a neighborhood at $x_0,$ satisfying that for each $x\in U_0, \, R(T_x)\cap N(T^+_0)=\{0\}$. We claim that there exists a neighborhood at $x_0,$ such that for each $x$ in it, $T_x$ has a generalized inverse $T^+_x$ satisfying $\lim\limits_{x\rightarrow x_0}T^+_x=T^+_0.$ Since $T_x$ from $X\rightarrow B(E,F)$ is continuous at $x_0$, we can assume $U_0\subset\{x\in X:\|T_x-T_0\|<\|T^+_0\|^{-1}\}$. Then the equivalence of the conditions $(i)$ and $(ii)$ in Theorem $1.1$ shows that for each $x$ in $U_0,$
$T^+ _x=T^+ _0C^{-1}_{T_0} (T ^ + _0, T_x)$ is a generalized inverse of $T_x,$ and $\lim\limits_{x\rightarrow
x_0}T^+_x=T^+_0$. Conversely, assume that for any generalized inverse $T_0^+$ of $T_0,$ there exists a neighborhood $U_0$ at $x_0,$ such that for each $x \in U_0,  T_x$ has a generalized inverse $T^+_x$ that satisfies $\lim\limits_{x\rightarrow x_0}T^+_x=T_0^+$.
 We claim that there exists a neighborhood $U$ at $x_0,$ such that $R(T_x)\cap N(T_0^+)=\{0\} \, \, \forall x\in U.$ Consider the operator valued map as follows,
$$P_x=I_E-T^+_x T_x\quad \mbox{for}\quad x\in U_0.$$
 Write $P_0=P_{x_0}.$ Obviously $P_0=I_E-T^+_0T_0, R(P_x)=N(T_x)$ and
 $R(P_0)=N(T_0).$ Write $V_0=\{x\in U_0:\|P_x-P_0\|<1\}\cap\{x\in
 U_0:\|T_x-T_0\|<\|T^+_0\|^{-1}\}.$
The question $4.11$ in [Ka] indicates that when $ \parallel P_x-P_0\parallel < 1,$  $$P_0 R(P_x)=R(P_0), \, \, \mbox{i.e.,} \, \, (I_E-T^+_0T_0)N(T_x)=N (T_0 ) \, \, \mbox{for} \, \, x \in V_0. $$ Then by the equivalence of the
 conditions $(i)$ and $(v)$ in Theorem $1.1$,
 $$R(T_x)\cap N(T^+_0)=\{0\},\quad\quad\forall x\in V_0.$$ The proof ends.

\textbf{$5$\quad The Proof of Theorem $1.5$}

Assume that subspaces $E_0$ and $E_1$ in $E$ possess a common complement $E_*$.
 First, go to verify that there is unique operator $\alpha $ in $B(E_0, E_*)$, satisfying $E_1=\{e+\alpha
e: \, \forall e\in E_0\}$.

Obversely, $$P^{E_*}_{E_0}P^{E_*}_{E_1}e=P^{E_*}_{E_0}(P^{E_*}_{E_1}e+P^{E_1}_{E_*})e=P^{E_*}_{E_0}e=e,\quad\forall
e\in E_0, \quad \eqno (5) $$ and
$$P^{E_*}_{E_1}P^{E_*}_{E_0}e=P^{E_*}_{E_1}(P^{E_*}_{E_0}e+P^ {E_0}_{E_*}e)=P^{E_*}_{E_1}e=e,\quad\forall
e\in E_1 .\quad \eqno (6)$$ The equalities $(5)$ and $(6)$ lead to $\alpha=\left.P^{E_0}_{E_*}P^{E_*}_{E_1}\right|_{E_0}$ fulfilling $E_1=\{e+\alpha e:\forall e\in E_0\}$. Indeed,  by $(6)$
$$e=P^{E_*}_{E_0}e+P^{E_0}_{E_*}e=P^{E_*}_{E_0}e+ \alpha P^{E_*}_{E_0}e \, \, \forall e\in E_1, \quad \eqno (7)$$ this indicates that
$E_1 \subset \{e + \alpha e: \, \, \mbox{for any } \, \, e \in E_0 \};$
by $(5)$ ,
$$e+\alpha e=P^{E_*}_{E_0}(P^{E_*}_{E_1}e)+\alpha P^{E_*}_{E_0}(P^{E_*}_{E_1}e) \, \, \forall e\in E_0;$$  by $(7),$
 $\{e + \alpha e: \, \, \mbox{for any } \, \, e \in E_0 \}\subset E_1 .$ Therefore, $E_1=\{e+\alpha
e: \, \forall e\in E_0\}$.

Assume that  $\beta$ is another operator in $B(E_0, E_*),$ satisfying that $E_1 = \{e + \beta e: \, \, \mbox{for any } \, \, e \in E_0\}. $  For any $h\in E_1$ let $h = e_0 + \alpha e_0=e_1 + \beta e_1$ where both $e_0$ and $e_1$ belong to $E_0.$ Then, $e_0=e_1$ and $\alpha e_0= \beta e_1$. This shows $\alpha = \beta,$ which means that $\alpha$ is unique. For any $\alpha \in B(E_0, E_*)$ let $E^\sharp = \{e + \alpha e :  \forall e \in E_0\}.$ We want to prove that $E=E^\sharp \oplus E_*,$ which includes the following three items: $(i)$  $E^\sharp $ is a closed subspace in $E,$ $(ii)$ $E^\sharp \cap E_*=\{0\},$ and $(iii)$ $E=E^\sharp + E_*.$

$(i)$.  Let $e_n+\alpha e_n\rightarrow e_*$ as
$n\rightarrow\infty$ where $e_n\in E_0, n=1,2,\cdots.$ Because the three operators $\alpha,$
$P^{E_*}_{E_0}$ and $P_{E_*}^{E_0}$ are all bounded linear operators,
$$\lim \limits_{n\rightarrow \infty}e_n=\lim \limits _{n\rightarrow \infty}P^{E_*}_{E_0}(e_n+\alpha e_n)= P^{E_*}_{E_0}e_*\in E_0, $$ and
$$\lim \limits _{n\rightarrow \infty}\alpha e_n=\lim \limits _{n\rightarrow \infty} P^{E_0}_{E_*}(e_n+\alpha e_n)
=\alpha P^{E_*}_{E_0}e_*.$$ Therefore,
$$e_*=\lim \limits_{n\rightarrow \infty}(e_n+ \alpha e_n)=P^{E_*}_{E_0}e_*+\alpha P^{E_*}_{E_0}e_*\in
E_0.$$ This shows that $E^\sharp $ is closed.

$(ii).$  For each $e \in E^\sharp \cap E_*$ let $e_0 $ in $E_0$ satisfy that $ e=e_0 + \alpha e_0.$  Then $e_0=0$ because $e\in E_*,$ and hence $e=0.$ This indicates that $E^\sharp\cap E_*=\{0\}.$

 $(iii).$ Obviously, in order to prove $(iii),$ it only needs to verify that $E^\sharp+E_*\supset E.$ Evidently, for each  $e\in E,$
$$e=P^{E_*}_{E_0}e+P^{E_0}_{E_*}e=(P^{E_*}_{E_0}e+\alpha P^{E_*}_{E_0}e)+(P^{E_0}_{E_*}e-\alpha P^{E_*}_{E_0}e),$$ where $P^{E_*}_{E_0}e+\alpha P^{E_*}_{E_0}e$ and $P^{E_0}_{E_*}e-\alpha P^{E_*}_{E_0}e$ belong to $E^\sharp $ and $E_*$ respectively. Therefore, $E^\sharp+ E_*\supset E$.
Theorem $1.5$ is proved.

\textbf{$6$ \quad Generalized regular point}

 Let $f$ be a $c^1$ map
from an open set $U$ in a Banach space $E$ to another Banach space $F$. It is well known that when the point $x_0$ is the submersion, immersion, and subimmersion points of $f$, $f$ has the submersion, immersion and subimmersion theorems at $x_0$, respectively. these three theorems provide an important way to study the local behavior of $f$ near the point $ x_0$ with using $f'(x_0)$. (Refer [Abr] and [Zei] .) In [Ber], Berger shows that it is not yet known whether the rank theorem in advanced calculus holds even if $f$ is Fredhlom map. For abbreviation, write $T_0=f'(x_0)$, and let $T_0^+$ denote a generalized inverse of $T_0$. In 1999, we proposed the following diffeomorphsms $u$ in $E$ and $v$ in $F$: $$ u(x)=T_0^+(f(x)-f(x_0))+(I_E-T_0^+T_0)(x-x_0)$$ and  $$ v(y)=(f\circ u^{-1}\circ T_0^+)(y)+(I_F-T_0T_0^+)y, $$ which satisfy that $u(x_0)=0, u'(x_0)=I_E, v(0)=f(x_0),$ and $ v'(0)=I_F.$ This proves that the formula
 $$f(x) =(v\circ f(x_0) \circ u)(x) \, \, \mbox{for all} \, \, x \in U_0$$ holds if and only if $x_0 $ is the generalized regular point of $f.$  In other words, the necessary and sufficient condition for $f$ to be locally conjugate to $f'(x_0)$ near $x_0$ is that $x_0$ is a generalized regular point of $f.$ It is called the complete rank theorem.(Refer to [Ma$1$ ], [Ma$8$], [Ber], and [Zei].) The generalized regular point produces many important results, such as the generalized pre-image theorem (see [Ma3]), extremum principle of real functional under the generalized regular constraint (see [Ma7]), the theorems $3.1 $ and $3.2$ of this paper, the complete rank theorem mentioned above, etc. Therefore, ones should recognize that the generalized regular point is a good mathematical concept after the regular point.

\textbf{Acknowlegment} I am fortunate to have completed this paper in my nineties and dedicate it to Tseng Yuanrong Functional Research Center and My Alma mater Anqing No.1 Middle School as a memorial.

\vskip 0.1cm
\begin{center}{\bf References}
\end{center}
\vskip -0.1cm
\medskip
{\footnotesize
\def\REF#1{\par\hangindent\parindent\indent\llap{#1\enspace}\ignorespaces}

\REF{[Abr]}\ R. Abraham, J. E. Marsden, and T. Ratin, Manifolds,
Tensor Analysis and Applications, 2nd ed., Applied Mathematical
Sciences 75, Springer, New York, 1988.

\REF{[An]}\ V. I. Arnol'd, Geometrical Methods in the Theory of
Ordinary Differential Equations, 2nd ed., Grundlehren der
Mathematischen Wissenschaften 250, Springer, New York, 1988.

\REF{[Ber]}\ M, Brger, Nonlinearity and Functional Analysis, New York:
Academic Press, 1976.

\REF{[Bo]}B.Boss, D.D.Bleecker, Topology and Analysis, The
Atiyah-Singer Index Formula and Gauge-Theoretic Physics, New York
Springer-Verlag, 1985.

\REF{[Caf]}\ V. Cafagra, Global invertibility and finite
solvability, pp. 1-30 in Nonlinear Functional Analysis (Nework, NJ,
1987), edited by P. S. Milojevic, Lecture Notes in Pure and Appl.
Math. 121, Dekker, New York, 1990.

 \REF{[Ka]}\ T. Kato,  Perturbation Theory for Linear
Operators, New York: Springer-Verlag, 1982.

\REF{[Ma1]}\ Jipu Ma, (1,2) inverses of operators between Banach
spaces and local conjugacy theorem, Chinese Ann. Math. Ser. B.
20:1(1999), 57-62.

\REF{[Ma2]}\ Jipu Ma, Rank theorem of operators between Banach
spaces, Sci China Ser. A 43(2000), 1-5.

\REF{[Ma3]}\ Jipu Ma, A generalized preimage theorem in global
analysis, Sci, China Ser. A 44:33(2001), 299-303.

\REF{[Ma4]}\ Jipu Ma, A rank theorem of operators between Banach
spaces, Front. of Math. in China 1(2006), 138-143.

\REF{[Ma5]}\ Jipu Ma, A generalized transversality in global
analysis, Pacif. J. Math.,  236:2(2008),  357-371.

\REF{[Ma6]}\ Jipu Ma, A geometry characteristic of Banach spaces
with $C^1$-norm, Front. Math. China,  2014, 9(5): 1089-1103.

\REF{[Ma7]} Jipu Ma, A principle for critical point under
generalized regular constrain and Ill-posed Lagrange multipliers,
Numerical Functional and Optimization Theory 36.225:(2015), 369.

\REF{[Ma8]}\  Jipu Ma, Local conjugacy  theorem, rank theorem in
advanced calculus and a generalized principle for constructing Banach
manifold, Sci. China Ser. A 43 (2000), 1233-1237.

\REF{[N-C]}\ M. Z. Nashed, X. Chen, Convergence of Newton-like methods for singular equations using outer
 inverses, Numer. Math., 66:(1993), 235-257.

\REF{[P]}\ R. Penrose, A generalized inverse for Matrices, Proc. Cambridge
Philos. Soc. 1955, 51:406-413.

\REF{[Zei]}\ A. E. Zeilder, Nonlinear Functional Analysis and its Applications,
IV. Applications to Mathematical Physics. New York: Springer -verlag, 1988.

%
%

\end{document}